\documentclass[12pt]{article}
\usepackage{amsmath,amsfonts,mathrsfs,wasysym,times,lineno}
\usepackage{amsthm,amscd,amsfonts,newlfont}
\usepackage{amssymb, upref, color}
\usepackage[dvips]{graphicx}
\usepackage{epsf, subfigure, verbatim}
\usepackage{latexsym}
\usepackage[colorlinks]{hyperref}
\usepackage{mathrsfs}
\usepackage{epsfig}
\usepackage{texdraw}
\usepackage[dvips]{psfrag}

\font\tenmsb=msbm10 \font\sevenmsb=msbm7
  \font\fivemsb=msbm5
  \catcode`\@=11
  \ifx\amstexloaded@\relax\catcode`\@=\active
  \endinput\else\let\amstexloaded@\relax\fi
  \def\spaces@{\space\space\space\space\space}
  \def\spaces@@{\spaces@\spaces@\spaces@\spaces@\spaces@}
  \def\space@.{\futurelet\space@\relax}
  \space@.%
  \def\Err@#1{\errhelp\defaulthelp@\errmessage{AmS-TeX error:#1}}
  \def\relaxnext@{\let\next\relax}
  \def\accentfam@{7}
  \def\noaccents@{\def\accentfam@{0}}
  \def\Cal{\relaxnext@\ifmmode\let\next\Cal@\else
  \def\next{\Err@{Use \string\Cal\space only in math mode}}\fi\next}
  \def\Cal@#1{{\Cal@@{#1}}}
  \def\Cal@@#1{\noaccents@\fam\tw@#1}
  \def\Bbb{\relaxnext@\ifmmode\let\next\Bbb@\else
  \def\next{\Err@{Use \string\Bbb\space only in math mode}}\fi\next}
  \def\Bbb@#1{{\Bbb@@{#1}}}
  \def\Bbb@@#1{\noaccents@\fam\msbfam#1}
  \newfam\msbfam
  \textfont\msbfam=\tenmsb
\scriptfont\msbfam=\sevenmsb
\scriptscriptfont\msbfam=\fivemsb

\@addtoreset{equation}{section}

\begin{document}

\title{{\bf ON THE  CUBIC PERTURBATIONS OF THE SYMMETRIC 8-LOOP HAMILTONIAN}\thanks {This work was partially supported by NNSF of China grant number 11771282 (C. Li) and by NNSF of China grant numbers 11431008 and 11771282 (J. Yu).}}
\author{{Iliya D. Iliev$^{a}$, \, Chengzhi Li$^{b}$ and Jiang Yu$^{c}$\,{\footnote{Corresponding Author. {\newline E-mail addresses: iliya@math.bas.bg, licz@math.pku.edu.cn, jiangyu@sjtu.edu.cn}}}}\\[2mm]
{\small\it $^{a}$Institute of Mathematics, Bulgarian Academy of Sciences, Bl. 8, 1113 Sofia, Bulgaria}\\[2mm]
{\small\it $^{b}$School of Mathematical Sciences, Peking University, Beijing 100871, PR China}\\[2mm]
{\small\it $^{c}$School of Mathematical Sciences, Shanghai Jiao Tong University, Shanghai 200240, PR China}}
\date{}
\maketitle{}

\begin{abstract} We study arbitrary cubic perturbations of the symmetric 8-loop Hamiltonian, which
are linear with respect to the small parameter. It is shown that when the first 4 coefficients in the
expansion of the displacement functions corresponding to both period annuli inside the loop vanish,
the system becomes integrable with the following three strata in the center manifold: Hamiltonian,
reversible in $y$ and reversible in $x$. In the latter case, the first integral is of Darboux type and
we calculate it explicitly.

Next we prove that the cyclicity of each of period annuli inside the loop is five and the total cyclicity
of both is at most nine. For this, we use Abelian integrals method together with careful study the geometry
of the separatrix solutions of related Riccati equations in connection to the second-order Melnikov functions.
\bigskip


\medskip


\end{abstract}

\section  {Introduction. }

We consider the symmetric eight-loop Hamiltonian
\begin{equation}\label{E1}
H(x,y)=\frac12y^2-\frac12x^2+\frac14 x^4,
\end{equation}
having two critical values: $h_0=0$, corresponding to the saddle $S$ at the origin, and  $h_1=-\frac14$,
corresponding to the centers at $C=(1, 0)$ and $C^*=(-1,0)$. We are going to study
the bifurcation of limit cycles from the two continuous families of ovals
$\mathcal{A}=\{\delta(h)\}$, $\mathcal{A}^*=\{\delta^*(h)\}\subset \{H=h\}$,  $h\in (-\frac14,0)$,
which surround respectively the centers $C$ and $C^*$, under general small cubic perturbations
\begin{equation}\label{E2}
\begin{array}{l} \dot{x}=H_y+\varepsilon f(x,y),\\
\dot{y}=-H_x+\varepsilon g(x,y)
\end{array}
\end{equation}
which are linear with respect to the small parameter $\varepsilon$. We assume that
$$f(x,y)=\sum_{i+j=0}^3 a_{ij} x^iy^j, \quad  g(x,y)=\sum_{i+j=0}^3 b_{ij} x^iy^j.$$
Our first goal is to calculate all the Melnikov functions $M_k(h)$ and  $M^*_k(h)$ in the expansion of the respective
displacement maps
\begin{equation}\label{E3}
d(h,\varepsilon)=\varepsilon M_1(h)+\varepsilon^2 M_2(h)+\varepsilon^3 M_3(h)+\ldots,
\end{equation}
$$d^*(h,\varepsilon)=\varepsilon M^*_1(h)+\varepsilon^2 M^*_2(h)+\varepsilon^3 M^*_3(h)+\ldots. \eqno(1.3^*)$$
which we need in order to determine the cyclicity of both annuli under the perturbation (\ref{E2}). For this,
we use the procedure proposed by Fran\c{c}oise \cite{fra} as it is most adequate to our purposes.
It turns out that $M_2(h)$ and $M_2^*(h)$ belong to most broad linear space among all $M_k(h)$ and $M_k^*(h)$.
Next, we wish to determine how many coefficients in (\ref{E3}) or $(1.3^*)$ should vanish in $(-\frac14, 0)$ in order to obtain
an integrable perturbation (that is, in order system (\ref{E2}) to have a center near $C$ or $C^*$, respectively).
And at last, taking a perturbation with $M_1(h)=M^*_1(h)=0$ and investigating the coefficients
$M_2(h)$, $M^*_2(h)$ in (\ref{E3}), $(1.3^*)$ that allow maximal number of zeros, we
are going to determine the cyclicity of the period annuli $\mathcal{A}$, $\mathcal{A}^*$ of $H$ as well as their total
cyclicity, that is the maximal number of limit cycles which can be produced from  $\mathcal{A}$, $\mathcal{A}^*$
and $\mathcal{A}\cup\mathcal{A}^*$ under the small perturbation (\ref{E2}).

There is a lot of papers considering small perturbations like (\ref{E2}), even in the case of arbitrary degree polynomials $f,g$.
However, most of them either investigate only the first Melnikov function $M_1(h)$ \cite{ChL}, or deal with some particular
low-degree perturbations aiming to obtain exact results, or both \cite{DLi}. Very few papers investigate polynomial
perturbations $f(x,y,\varepsilon)$, $g(x,y,\varepsilon)$ which are analytic with respect to the small parameter $\varepsilon$.
For some results concerning specific low-degree perturbations of this kind see \cite{ip}, \cite{HamGarGav} and for results in
a general setting see \cite {GI}.

If $h\in\Sigma=(-\frac14,0)$, we denote $I_k(h)=\int_{\delta(h)}x^kydx$, $k=0,1,\ldots$, assuming that the integrals have
the standard positive (counterclockwise) orientation. In the statements that follow, the coefficients of $M_k(h)$, $M^*_k(h)$
(denoted by Greek letters) are different for different $k$ and all they are explicitly determined by the coefficients
$a_{ij}, b_{ij}$ of the perturbation  (\ref{E2}).  For more details see Propositions 1--4 in the next sections, which together
imply the proof of Theorem 1.  Our main results in this paper are the following.

\vspace{2ex}
\noindent
{\bf Theorem 1.} {\it The following statements hold about $(\ref{E3})$ and $(1.3^*)$:}

\vspace{1ex}
\noindent
(i)  {\it The first coefficients $M_1(h)$ and  $M^*_1(h)$ take the form:}
$$\begin{array}{l}
 M_1(h)=(\alpha_0+\alpha_1h)I_0(h)+\beta_0I_1(h)+\gamma_0I_2(h),\quad h\in \Sigma,\\[2mm]
M^*_1(h)=(\alpha_0+\alpha_1h)I_0(h)-\beta_0I_1(h)+\gamma_0I_2(h),\quad h\in \Sigma,
\end{array}$$

\noindent
(ii) {\it If the first coefficients in $(\ref{E3})$ and $(1.3^*)$ vanish, then the second ones take the form:}
$$\begin{array}{l}
 M_2(h)=(\alpha_0+\alpha_1h)I_0(h)+(\beta_0+\beta_1h)I_1(h)+(\gamma_0+\gamma_1h)I_2(h),\quad h\in \Sigma,\\[2mm]
M^*_2(h)=(\alpha_0+\alpha_1h)I_0(h)-(\beta_0+\beta_1h)I_1(h)+(\gamma_0+\gamma_1h)I_2(h),\quad h\in \Sigma,
\end{array}$$

\noindent
(iii) {\it If the first two coefficients in $(\ref{E3})$ and $(1.3^*)$ vanish, then the third ones take the form:}
$$\begin{array}{l}
 M_3(h)=(\alpha_0+\alpha_1h)I_0(h)+(\beta_0+\beta_1h)I_1(h)+\gamma_0I_2(h),\quad h\in \Sigma,\\[2mm]
M^*_3(h)=(\alpha_0+\alpha_1h)I_0(h)-(\beta_0+\beta_1h)I_1(h)+\gamma_0I_2(h),\quad h\in \Sigma,
\end{array}$$

\noindent
(iv) {\it If the first three coefficients in $(\ref{E3})$ and $(1.3^*)$ vanish, then the fourth coefficients
$M_4(h)$ and  $M^*_4(h)$, if not zero, form three one-dimensional subspaces spanned over $\{I_1, hI_1, I_2\}$.}

\vspace{2ex}

\noindent
(v) {\it If the first four coefficients in $(\ref{E3})$ and $(1.3^*)$ vanish, then system $(\ref{E2})$ is integrable
and belongs to one of the three strata: Hamiltonian, reversible in $y$, reversible in $x$.}

\vspace{2ex}
Then investigation the number of zeros of $M_2(h)$ and $M_2^*(h)$ in $\Sigma$  yields:

\vspace{2ex}
\noindent
{\bf Theorem 2.} {\it  For any compact region $K$ contained in $\mathcal{A}$ or $\mathcal{A^*}$,
system $(\ref{E2})$ has for small $\varepsilon$ at most 5 limit cycles in $K$, including their multiplicities.
This bound is exact.}

\vspace{2ex}
\noindent
What concerns the simultaneous bifurcation of limit cycles from both annuli, we prove:

\vspace{2ex}
\noindent
{\bf Theorem 3.} {\it System $(\ref{E2})$ can have}

\vspace{1ex}
\noindent
(i) {\it Up to 9 limit cycles, with possible maximal distributions 5+4; 4+5 if all 6 parameters in $M_2(h)$ are present.}

\vspace{1ex}
\noindent
(ii) {\it Up to 7 limit cycles, with possible  maximal distributions 3+4; 4+3 if  $\beta_1=0$ in $M_2(h)$.}

\vspace{1ex}
\noindent
(iii) {\it Up to 6 limit cycles, with possible maximal distribution 3+3  if  $\beta_1=\beta_0=0$ in $M_2(h)$.}

\vspace{2ex}
By using the explicit expansions of the integrals $I_k(h)$ at level $h=-\frac14$ corresponding to the centers, we obtain in addition
the following estimates for the total cyclicity of $C$ and $C^*$ with respect to perturbation (1.2), which also yield some estimates from
below about the number of limit cycles in the system.

\vspace{2ex}
\noindent
{\bf Theorem 4.} {\it The following maximal distributions of small-amplitude limit cycles in $(1.2)$ are possible: 5+0, 0+5, 4+1, 1+4, 3+3.}

\vspace{2ex}
We would like to note that when calculated for degree $n$ perturbations, the function $M_1(h)$ in (3) has coefficients
$p_k(h)$ at $I_k$, $k=0,1,2$ which are also polynomials in $h$. However, the degree of $p_0$ is always strictly greater
than the degree of $p_2$ \cite{Petrov}, apart of the situation with $M_2(h)$ here. Hence, no known estimate about
$M_1(h)$ can be used to obtain the {\it exact} bounds for $\mathcal{A}$ and $\mathcal{A}^*$, stated in Theorems 2 and 3.

Our proofs are purely classical and use simple tools only, such as two-dimensional Fuchsian systems and the respective
Riccati equations, the geometric properties of the separatrix curves in their phase plane, contact points of the flow with
some simple curves, etc.


\section  {Calculation of the coefficients $M_1(h)$ and $M^*_1(h)$. }

If $h\in\Sigma=(-\frac14,0)$, we denote $I_k(h)=\int_{\delta(h)}x^kydx$, $k=0,1,\ldots$ and
 $I_{kl}(h)=\int_{\delta(h)}x^ky^ldx$, $k,l=0,1,\ldots$. We will assume that the line integrals are oriented  in a
counterclockwise direction.
It is well known that for odd $l$, the integrals $I_{kl}(h)$ are expressed as polynomial (in $h$) envelopes of the three
basic integrals $I_0(h)$, $I_1(h)$ and $I_2(h)$, while for $l$ even, they do vanish because of the symmetry of the ovals
$\delta(h)$ with respect to $y$. It is well known that if one writes the system (\ref{E2}) in a Pfaffian form
$$dH=\varepsilon\omega,\quad \omega=g(x,y)dx-f(x,y)dy,$$
then $M_1(h)$ is given by
\begin{equation}\label{E4}
M_1(h)=\int_{\delta(h)}\omega,\quad h\in\Sigma.
\end{equation}
Moreover, one can verify that the so called {\it star property} holds in our case. Namely, if $M_1(h)$ vanishes identically
 in $\Sigma$, then $\omega$ takes the form $\omega = dQ(x,y)+q(x,y)dH$ with appropriate {\it polynomials} $Q,q$
(which are not uniquely defined). If so, then the next coefficient $M_2(h)$ is simply
\begin{equation}\label{E5}
M_2(h)=\int_{\delta(h)}q(x,y)\omega,\quad h\in\Sigma
\end{equation}
and so on. For more details about star property and the recursive procedure, see  \cite{fra}.

Of course, this procedure could be generalized for the cases which do not possess the star property, but then the price
to pay is that $Q$ and $q$ will no more be polynomials. Nevertheless $M_k(h)$ can still be Abelian integrals \cite{GI}.
By the way, the third period annulus of $H$: $\{\delta(h)\}\subset \{H=h\}$, $h>0$ which surrounds the eight loop
$\{H=0\}$ does not obey the star property. This fact affects the expressions for the higher order coefficients $M_k(h)$,
which for $h>0$ are {\it rational} envelopes of the two integrals  $I_0(h)$, $I_2(h)$ only, since $I_1(h)$ becomes zero
here.

In other words, the star property should in general be attributed to a couple $(H,\mathcal{A})$ and not to the Hamiltonian
function $H$ alone.

Now we are ready to formulate

\vspace{2ex}
\noindent
{\bf Proposition 1.} {\it The following statements hold about $(\ref{E3})$ and $(1.3^*)$:}

\vspace{1ex}
\noindent
(i)  {\it The first coefficients $M_1(h)$ and  $M^*_1(h)$ take the following form:
$$\begin{array}{l}
 M_1(h)=(\alpha_0+\alpha_1h)I_0(h)+\beta_0I_1(h)+\gamma_0I_2(h),\quad h\in \Sigma,\\[2mm]
M^*_1(h)=(\alpha_0+\alpha_1h)I_0(h)-\beta_0I_1(h)+\gamma_0I_2(h),\quad h\in \Sigma,
\end{array}$$
where}
$\alpha_0=a_{10}+b_{01}$, $\;\alpha_1=\frac47(a_{12}+3b_{03})$, $\;\beta_0=2a_{20}+b_{11}$,
$\gamma_0=3a_{30}+b_{21}+\frac17(a_{12}+3b_{03}).$

\vspace{1ex}
\noindent
(ii)  {\it One has} $M_1(h)\equiv 0 \Leftrightarrow M^*_1(h)\equiv 0\Leftrightarrow
   \alpha_0=\alpha_1=\beta_0=\gamma_0=0 \Leftrightarrow$
$a_{10}+b_{01}=a_{12}+3b_{03}=2a_{20}+b_{11}=3a_{30}+b_{21}=0.$

\vspace{2ex}
\noindent
{\bf Proof.} (i) By (\ref{E4}), we have  $M_1(h)=-\int\!\!\int_{H<h}(f_x+g_y)dxdy=\int_{\delta(h)}G(x,y)dx$ where
$$\textstyle{ G=\alpha_0y+\beta_0xy+(\frac12a_{11}+b_{02})y^2
+(3a_{30}+b_{21})x^2y+(a_{21}+b_{12})xy^2+(\frac13a_{12}+b_{03})y^3}$$
Using identity (I-1) from the Appendix, we obtain the needed formula. The expression for $M^*_1(h)$ is
derived by changing the variables $(x,y)\to (-x,-y)$ in (\ref{E2}) which moves $C^*$ to $(1,0)$.
This transformation only changes signs of the coefficients  $a_{ij}$, $b_{ij}$ when $i+j$ is even. Therefore it
affects the sign of $\beta_0$ only. Statement (ii) is a consequence of the fact that in $\Sigma$, each $\mathbb{R}[h]$
module over the basic integrals $I_0, I_1, I_2$  is free. $\Box$


\section  {Calculation of the coefficients $M_2(h)$ and $M^*_2(h)$. }

Assuming that $M_1(h)$ vanishes and using Proposition 1 (ii), we can perform elementary calculations to rewrite
$\omega$ in the form
\begin{equation}\label{E6}
\omega=dQ(x,y)+(\lambda y^2+\mu xy^2)dx
\end{equation}
where we have denoted for short
$$\textstyle \lambda=\frac12a_{11}+b_{02}, \quad\mu=a_{21}+b_{12},\quad \mbox{\rm and}\quad  Q=Q_1+Q_2,$$
$$\begin{array}{rl}
Q_1 &=-(a_{00}y+a_{10}xy+a_{20}x^2y+a_{30}x^3y+\frac13a_{02}y^3+\frac13a_{12}xy^3),\\[2mm]
Q_2 &=b_{00}x+\frac12b_{10}x^2+\frac13b_{20}x^3+\frac14b_{30}x^4-\frac12a_{01}y^2
-\frac12a_{11}xy^2-\frac12a_{21}x^2y^2-\frac14a_{03}y^4.\end{array}$$
If it happens that both $\lambda,\mu$ are zero, we have nothing more to do since system (\ref{E2}) becomes Hamiltonian
in this case, with a Hamiltonian function $H(\varepsilon)=H-\varepsilon Q$.

Thus we will assume in what follows that $|\lambda|+|\mu|\neq 0$. If so, by (\ref{E6}) and formulas (II-1) and (II-2)
from the Appendix, we conclude that one can apply (\ref{E5}), with $q=q_1=-(2\lambda x+\mu x^2)$.
Therefore,
$$M_2(h)=\int_{\delta(h)}q_1\omega=2\int_{\delta(h)}(\lambda+\mu x)Q_1dx$$
since all other integrals do vanish by symmetry in $y$. This yields
$$\begin{array}{rl}
M_2(h) &=-2[a_{00}\lambda I_0+(a_{10}\lambda+a_{00}\mu)I_1+(a_{20}\lambda+a_{10}\mu)I_2
+(a_{30}\lambda+a_{20}\mu)I_3\\[2mm]
& +a_{30}\mu I_4+\frac13a_{02}\lambda I_{03}+(\frac13a_{12}\lambda+\frac13a_{02}\mu)I_{13}
+\frac13a_{12}\mu I_{23}]\end{array}$$

\vspace{2ex}
\noindent
{\bf Proposition 2.} {\it Assume that the first coefficients in $(\ref{E3})$ and $(1.3^*)$ vanish. Then:}

\vspace{2ex}
\noindent
(i)  {\it The second coefficients $M_2(h)$ and  $M^*_2(h)$ take the following form:}
$$\begin{array}{l}
 M_2(h)=(\alpha_0+\alpha_1h)I_0(h)+(\beta_0+\beta_1h)I_1(h)+(\gamma_0+\gamma_1h)I_2(h),\quad h\in \Sigma,\\[2mm]

M^*_2(h)=(\alpha_0+\alpha_1h)I_0(h)-(\beta_0+\beta_1h)I_1(h)+(\gamma_0+\gamma_1h)I_2(h),\quad h\in \Sigma,
\end{array}$$

\noindent
{\it where the six Greek letter coefficients are independently free and given by}

\vspace{2mm}
$\alpha_0=-2a_{00}\lambda, \quad \alpha_1=-\frac87a_{02}\lambda-(\frac87a_{30}+\frac{8}{63}a_{12})\mu$,

\vspace{1mm}
$\beta_0=-2(a_{10}+a_{30}+\frac18a_{12})\lambda-2(a_{00}+a_{20}+\frac18a_{02})\mu, \quad
\beta_1=-a_{12}\lambda-a_{02}\mu$,

\vspace{1mm}
$\gamma_0=-2(a_{20}+\frac17a_{02})\lambda-2(a_{10}+\frac87a_{30}+\frac{8}{63}a_{12})\mu, \quad
\gamma_1=-\frac89a_{12}\mu.$

\vspace{2mm}
\noindent
(ii)  {\it The second coefficients $M_2(h)$ and  $M^*_2(h)$ vanish if and only if one of the five conditions holds:}

(a)  $Q_1=a_{00}(x^2-1)y,\qquad\lambda=0;$

(b)  $Q_1=a_{10}(x^3-x)y,\qquad\mu=0;$

(c)  $Q_1=a_{10}(x^2-x)y,\qquad\mu=\lambda;$

(d)  $Q_1=-a_{10}(x^2+x)y,\quad\mu=-\lambda;$

(e) $Q_1=0,\qquad \lambda\mu (\lambda^2-\mu^2)\neq 0.$

\vspace{2ex}
\noindent
{\bf Proof.} (i) We use formulas (I-1) and (II-3)-(II-6) in the Appendix to express the integrals $I_{kl}$ which
appear in $M_2(h)$ through $I_k$, $hI_k$, $k=0,1,2$. Thus we obtain easily the needed expression of $M_2(h)$.
Recalling that change the places of $C$ and $C^*$ results in changing signs of the coefficients  $a_{ij}$,
$b_{ij}$ when $i+j$ is even, and therefore of $\lambda$ and $\beta_0$, $\beta_1$ only, we obtain the
expression of $M_2^*$. (ii). The system $\alpha_i=\beta_i=\gamma_i=0$, $i=0,1$ is equivalent to
$$\begin{array}{r}a_{00}\lambda=0\\
a_{02}\lambda+a_{30}\mu=0\\
(a_{10}+a_{30})\lambda+(a_{00}+a_{20})\mu=0\\
a_{12}\lambda+a_{02}\mu=0\\
a_{20}\lambda+(a_{10}+a_{30})\mu=0\\
a_{12}\mu=0
\end{array}$$
Considering this as a system about $a_{jk}$, we note that its determinant vanishes if $\lambda\mu(\mu^2-\lambda^2)=0$
which yields (e). The solutions over the zero set are as shown in (a)-(d). $\Box$


\section  {Calculation of the coefficients $M_3(h)$ and $M^*_3(h)$. }

Let us begin by noticing that in case (e) the system (\ref{E2}) becomes reversible (that is symmetric with respect to $y$)
and therefore we have nothing more to do since all coefficients  $M_k(h)$ and $M^*_k(h)$ with indices $k\geq 3$ will be
zero as well. Indeed, in case (e), system (\ref{E2}) reduces to the equation
$$\frac12\xi'=\frac{x-x^3+\varepsilon(b_{00}+b_{10}x+b_{20}x^2+b_{30}x^3+(b_{02}+b_{12}x)\xi)}
{1+\varepsilon(a_{01}+a_{11}x+a_{21}x^2+a_{03}\xi)}$$
with respect to $\xi=y^2$. If $F_\pm(x,\xi)=h$ is the solution of the equation which exists for small $\varepsilon$
in a neighborhood of the point $(x,\xi)=(\pm 1,0)$, then $F_\pm(x,y^2)=h$ is a first integral of (\ref{E2}). By symmetry,
the two foci $(x_\pm, 0)$ near $(\pm 1,0)$ are centers.

For the other cases (a)-(d), we first need to calculate the respective $q_2$, a function such that $q_1\omega \sim q_2dH$
modulo exact forms. Then
$$M_3(h)=\int_{\delta(h)}q_2\omega.$$
We shall look for $q_2$ in the form $q_2=q_{20}+q_{21}+q_{22}$ where $q_{20}$ comes from the one-form (see (\ref{E6}))
$\omega_0=q_1(\lambda y^2+\mu xy^2)dx$, $q_{21}$ comes from $\omega_1=q_{1}dQ_1$ and  $q_{22}$ comes from
 $\omega_2=q_{1}dQ_2$. Thus,
$$\textstyle
\omega_0=-(2\lambda x+\mu x^2)(\lambda+\mu x)y^2dx\sim -2Hd(\lambda^2 x^2+\lambda\mu x^3+\frac14\mu^2x^4)\sim
\frac12 q_1^2dH=q_{20}dH$$
In order to consider all cases simultaneously, let us denote $m=2a_{00}\mu$ in case (a) and $m=2a_{10}\lambda$ in cases
(b), (c), (d). We will assume in what follows that $m\neq 0$. Elsewhere, one obtains that (\ref{E2}) is either Hamiltonian or
reversible, as we mentioned above. Then, an easy calculation using (a)-(d)  yields
$$\omega_1 \sim 2(\lambda+\mu x)Q_1dx  \sim m(x^3-x)ydx\sim mydH=q_{21}dH,$$
 the latter following from formula (II-3).
Finally,
$$\textstyle\omega_2\sim -(\lambda+\mu x)(a_{01}y^2+
a_{11}xy^2+a_{21}x^2y^2+\frac12a_{03}y^4)dx.$$
Using formulas (II-1), (II-2) and (III-1)-(III-4) from the Appendix, we immediately obtain
\begin{equation}\label{E7}
\begin{array}{rl}\omega_2 &\sim \lambda[2a_{01}x+a_{11}x^2+\frac23a_{21}x^3+\frac12a_{03}(8xH+\frac43x^3-\frac25x^5)]dH\\[2mm]
     & +\mu[a_{01}x^2+\frac23a_{11}x^3+\frac12a_{21}x^4+\frac12a_{03}(4x^2H+x^4-\frac13x^6)]dH=q_{22}dH.\end{array}
\end{equation}
Therefore
$$M_3(h)=\int_{\delta(h)}(q_{20}+q_{21}+q_{22})\omega=\int_{\delta(h)}(q_{20}+q_{22})d Q_1
+\int_{\delta(h)}q_{21}[dQ_2+(\lambda+\mu x)y^2dx]$$
$\hspace{12mm}=J_1+J_2$

\noindent
since all other integrals vanish by symmetry.

\vspace{2ex}
\noindent
{\bf Proposition 3.} {\it Assume that $m\neq 0$ and the first two coefficients in $(\ref{E3})$ and $(1.3^*)$ vanish. Then:}

\vspace{2ex}
\noindent
(i)  {\it The third coefficients $M_3(h)$ and  $M^*_3(h)$ take the following form:}
$$\begin{array}{l}
 M_3(h)=(\alpha_0+\alpha_1h)I_0(h)+(\beta_0+\beta_1h)I_1(h)+\gamma_0I_2(h),\quad h\in \Sigma,\\[2mm]

M^*_3(h)=(\alpha_0+\alpha_1h)I_0(h)-(\beta_0+\beta_1h)I_1(h)+\gamma_0I_2(h),\quad h\in \Sigma,
\end{array}$$

\noindent
{\it where the five Greek letter coefficients are independently free. Explicitly, }

$$\begin{array}{ll}
\alpha_0=mb_{00}, & \alpha_1=m(\frac47\lambda-\frac67a_{11}),\\[2mm]
\beta_0=m(b_{10}+b_{30}+\frac18\mu-\frac38a_{21}), &\beta_1=m(\frac12\mu-\frac32a_{21}),\\[2mm]
\gamma_0=m(b_{20}+\frac17\lambda-\frac{3}{14}a_{11}). &
\end{array}$$

\vspace{2mm}
\noindent
(ii)  {\it The third coefficients $M_3(h)$ and  $M^*_3(h)$ vanish if and only if
\begin{equation}\label{E8}
b_{00}=b_{20}=b_{10}+b_{30}=2\lambda-3a_{11}=\mu-3a_{21}=0,
\end{equation}
which is equivalent to}
\begin{equation}\label{E9}
\textstyle Q_2=\frac12(b_{10}-a_{01})y^2-b_{10}H-\frac 16(2\lambda x+\mu x^2)y^2-\frac14 a_{03}y^4.
\end{equation}

\vspace{2ex}
\noindent
{\bf Proof.} (i)  In order to handle all four cases (a)-(d) simultaneously, we keep both $\lambda$ and $\mu$.
Let us consider in detail first $J_1$. Clearly,
$$J_1=-\int_{\delta(h)}(q_{20}'+q_{22}')Q_1dx$$
where $'$ means a differentiation with respect to $x$ whilst $H=h$ is a constant over $\delta(h)$.
We easily obtain that
\begin{equation}\label{E10}
\begin{array}{l}q_{20}'+q_{22}' =\\[2mm]
2(\lambda+\mu x)[a_{01}+(a_{11}+2\lambda)x+(a_{21}+\mu)x^2+2a_{03}(H+\frac12x^2-\frac14x^4)],\\[2mm]
2(\lambda+\mu x)Q_1 =m(x^3-x)y.\end{array}
\end{equation}

Denote for a while by $K$ the expression in the brackets. Then
$$
\begin{array}{rl}J_1 &=-m\int_{\delta(h)}K(x^3-x)ydx=-m\int_{\delta(h)}Kyd(H-\frac12y^2)
=\frac13m\int_{\delta(h)}Kdy^3\\[2mm]
 &=-\frac13m\int_{\delta(h)}K'y^3dx=-\frac13m\int_{\delta(h)}[(a_{11}+2\lambda)+(2a_{21}+2\mu)x+2a_{03}(x-x^3)]y^3dx.
\end{array}
$$
Therefore
$$J_1=-\frac{m}{3}\int_{\delta(h)}[(a_{11}+2\lambda)+(2a_{21}+2\mu)x]y^3dx$$
as the integral at $a_{03}$ vanishes. Applying formulas (I-1) and (II-5) from the Appendix, we derive the formula
$$\textstyle J_1=-m[(a_{11}+2\lambda)(\frac17I_2+\frac47hI_0)+(a_{21}+\mu)(\frac14+h)I_1].$$

\vspace{2ex}
Consider now $J_2$.  We obtain, up to exact forms,

\vspace{2ex}
$q_{21}[dQ_2+(\lambda+\mu x)y^2dx]\sim$

\vspace{1ex}
$my(b_{00}+b_{10}x+b_{20}x^2+b_{30}x^3)dx-myd(\frac12a_{11}xy^2+\frac12a_{21}x^2y^2)+m(\lambda+\mu x)y^3dx$.

\vspace{1ex}
\noindent
The second form is equivalent to  $\frac16m(a_{11}x+a_{21}x^2)dy^3\sim -\frac16m(a_{11}+2a_{21}x)y^3dx.$

\vspace{1ex}
\noindent
Now applying formulas (I-1), (II-3), (II-5) from the Appendix and integrating the respective forms, one obtains
$$\begin{array}{rl}
J_2 &=m[b_{00}I_0+b_{10}I_1+b_{20}I_2+b_{30}I_3+(\lambda-\frac16a_{11})I_{03}+(\mu-\frac13a_{21})I_{13}]\\[2mm]
      &=m[b_{00}I_0+(b_{10}+b_{30})I_1+b_{20}I_2+(\lambda-\frac16a_{11})(\frac37I_2+\frac{12}{7}hI_0)
+(\mu-\frac13a_{21})(\frac38+\frac32h)I_1]\end{array}$$
Together, the expressions of $J_1$ and $J_2$ prove (i). The formula of $M_3^*$ follows by exchanging the centers.
The statement in (ii) is a consequence of (i) and the fact that all the coefficients should vanish. $\Box$


\section  {Calculation of the coefficients $M_4(h)$ and $M^*_4(h)$. }

Assume that  $M_k(h)$, $M^*_k(h)$ vanish for $k\leq 3$.
As in the previous section, to calculate $M_4(h)$ in the cases (a)-(d), we first need to find a function $q_3$ so that
 $q_2\omega \sim q_3dH$ modulo exact forms. Then
$$M_4(h)=\int_{\delta(h)}q_3\omega.$$

\vspace{2ex}
\noindent
{\bf Proposition 4.} {\it Assume that $m\neq 0$ and the first three coefficients in $(\ref{E3})$ and $(1.3^*)$ vanish. Then:}

\vspace{2ex}
\noindent
(i)  {\it The fourth coefficients $M_4(h)$ and  $M^*_4(h)$ take the form:}
$$\begin{array}{l}
 M_4(h)=M^*_4(h)=0   \quad  \mbox{\it if}\;\; \lambda=0,\\[2mm]
M_4(h)=6a_{10}^3\lambda(\frac14+h)I_1, \quad M^*_4(h)=-6a_{10}^3\lambda(\frac14+h)I_1 \quad  \mbox{\it if}\;\; \mu=0,\\[2mm]
M_4(h)=6a_{10}^3\lambda(	I_1-I_2), \quad M^*_4(h)=-6a_{10}^3\lambda(I_1+I_2) \quad  \mbox{\it if}\;\; \mu=\lambda,\\[2mm]
M_4(h)=6a_{10}^3\lambda(	I_1+I_2), \quad M^*_4(h)=-6a_{10}^3\lambda(I_1-I_2) \quad  \mbox{\it if}\;\; \mu=-\lambda.
\end{array}$$

\vspace{2mm}
\noindent
(ii)  {\it If $\lambda=0$, then system $(\ref{E2})$ becomes symmetric (reversible) with respect to $x$ and has centers at $(\pm 1,0)$
 for $\varepsilon$ sufficiently small.}

\vspace{2ex}
\noindent
{\bf Proof.} (i)
To begin with, we first recall that we have split $q_2$ into $q_2=q_{20}+q_{21}+q_{22}$.  According to this, we shall
split the one form we deal with as $q_2\omega=\omega_3+\omega_4$ where
$$\begin{array}{rl}
\omega_3 &=q_{21}[dQ_2+(\lambda y^2+\mu xy^2)dx]+(q_{20}+q_{22})dQ_1,\\[2mm]
\omega_4 &=q_{21}dQ_1+(q_{20}+q_{22})[dQ_2+(\lambda y^2+\mu xy^2)dx].\end{array}$$
We recall the reader that at this stage (see (\ref{E1}), (\ref{E6}), (a)-(d), (\ref{E7})-(\ref{E9})) we have
$$\begin{array}{l}
\omega=d(Q_1+Q_2)-\frac12y^2dq_1,\\[2mm]
q_1=-(2\lambda x+\mu x^2),\\[2mm]
Q_1=\bar{Q}_1y, \\[2mm]
Q_2=\frac12(b_{10}-a_{01})y^2-b_{10}H+\frac16q_1y^2-\frac14a_{03}y^4,\\[2mm]
q_{21}=my,\\[2mm]
q_{20}+q_{22}=\frac23q_1^2-(a_{01}+2a_{03}H)q_1+\varphi,\\[2mm]
\varphi=a_{03}(\frac23\lambda x^3+\frac12\mu x^4-\frac15\lambda x^5-\frac16\mu x^6),
\quad d\varphi=-2a_{03}(\frac12x^2-\frac14x^4)dq_1.
\end{array}$$
First, we easily obtain
$$\textstyle q_{21}[dQ_2+(\lambda y^2+\mu xy^2)dx]\sim \frac89m(\lambda+\mu x)y^3dx-mb_{10}ydH.$$
Here and below, the reader could verify relations between one-forms like $\omega_1\sim\omega_2$ which are not derived
 in full detail by expressing $\omega_1-\omega_2$ as $Fdx+Gdy$ and simply checking that $F_y=G_x$.
Next, we use the first equation in (\ref{E10}) to obtain
$$\begin{array}{l}
(q_{20}+q_{22})dQ_1\sim-Q_1d(q_{20}+q_{22})=-Q_1[(q_{20}'+q_{22}')dx+\partial_Hq_{22}dH]\\[2mm]
=-Q_1[2(\lambda+\mu x)Kdx+2a_{03}(2\lambda x+\mu x^2)dH]\end{array}$$
where $K$ is the function we used in (\ref{E10}) which by Proposition 3 (ii) is expressed now  as
$K=a_{01}+\frac83\lambda x+\frac43\mu x^2+a_{03}y^2$.
Then using the second equation in (\ref{E10}) we obtain
$$\begin{array}{l}
-2(\lambda+\mu x)Q_1 Kdx=-m(x^3-x)yKdx=-myKd(H-\frac12y^2)\\[2mm]
\sim -myKdH+m(\frac83\lambda x+\frac43\mu x^2)y^2dy\sim -myKdH-\frac89m(\lambda+\mu x)y^3dx.\end{array}$$
Summing up, we come to the relation $\omega_3\sim q_{31}dH$ where
\begin{equation}\label{E11}
\textstyle q_{31}=-m[(a_{01}+b_{10})y-\frac43q_1y+a_{03}y^3]+2a_{03}q_1Q_1.
\end{equation}
What about $\omega_4$, we similarly come to
$\omega_4\sim q_{30}dH$. In more detail,

$$\textstyle q_{21}dQ_1\sim -m\bar{Q}_1dH, $$

$$\begin{array}{l}(q_{20}+q_{22})[dQ_2+(\lambda y^2+\mu xy^2)dx]\sim -Q_2d(q_{20}+q_{22})-\frac12y^2(q_{20}+q_{22})dq_1\\[2mm]
=2a_{03}q_1Q_2dH+[Q_2(-\frac43q_1+a_{01}+a_{03}y^2)-y^2(\frac13q_1^2-\frac12a_{01}q_1-a_{03}Hq_1+\frac12\varphi)]dq_1\\[2mm]
=2a_{03}q_1Q_2dH+[-\frac14a_{03}^2y^6+(\frac12a_{03}q_1+\frac12a_{03}b_{10}-\frac34a_{01}a_{03})y^4]dq_1\\[2mm]
+[-\frac59q_1^2+(\frac43a_{01}-\frac23b_{10}+a_{03}H)q_1+\frac12a_{01}b_{10}-\frac12a_{01}^2-a_{03}b_{10}H-
\frac12\varphi]y^2dq_1\\[2mm]
+(\frac43b_{10}q_1-a_{01}b_{10})Hdq_1
\sim [-\frac14a_{03}^2y^6+(\frac12a_{03}q_1+\frac12a_{03}b_{10}-\frac34a_{01}a_{03})y^4]dq_1\\[2mm]
+[\frac{10}{27}q_1^3-(\frac43a_{01}+a_{03}H+\frac16a_{03}y^2)q_1^2+(a_{01}^2+(2a_{03}b_{10}-a_{01}a_{03})y^2-\frac12a_{03}^2y^4)q_1\\[2mm]
+\int\varphi q_1'dx+2a_{03}b_{10}\int q_1(x^3-x)dx-a_{03}\int q_1^2(x^3-x)dx]dH.\end{array}$$
Above, we have used relations
$$\textstyle -\frac12\varphi y^2dq_1=-\frac12\varphi q_1'y^2dx=-\frac12y^2d\int\varphi q_1'dx\sim (\int\varphi q_1'dx)dH$$
and similarly,
$$\begin{array}{l}a_{03}H(q_1-b_{10})y^2dq_1\sim-a_{03}(\frac12q_1^2-b_{10}q_1)dHy^2\\[1mm]
=-a_{03}(\frac12q_1^2-b_{10}q_1)y^2dH-a_{03}(q_1^2-2b_{10}q_1)Hd(H+\frac12x^2-\frac14x^4)\\[1mm]
\sim a_{03}(2b_{10}q_1-q_1^2)(\frac12y^2+H)dH+a_{03}[\int(2b_{10}q_1-q_1^2)(x^3-x)dx]dH.
\end{array}$$

\vspace{2ex}
And finally, we apply (IV-1)-(IV-3) to the remaining terms with $dq_1$ to obtain
\begin{equation}\label{E12}
\begin{array}{rl}
q_{30} &=
\frac{10}{27}q_1^3-(\frac43a_{01}+a_{03}H+\frac76a_{03}y^2)q_1^2+(a_{01}^2+2a_{01}a_{03}y^2
+a_{03}^2y^4)q_1\\[2mm]
& -m\bar{Q}_1+\int\varphi q_1'dx+(6a_{01}a_{03}-2a_{03}b_{10}+12a_{03}^2H)\int q_1(x^3-x)dx\\[2mm]
& -3a_{03}\int q_1^2(x^3-x)dx -12a_{03}^2\int q_1(x^3-x)(\frac14x^4-\frac12x^2)dx.
\end{array}
\end{equation}

\vspace{2ex}
By (\ref{E11}) and (\ref{E12}) therefore one obtains
$$M_4(h)=\int_{\delta(h)}(q_{30}+q_{31})\omega=\int_{\delta(h)}q_{30}dQ_1+\int_{\delta(h)}q_{31}
[dQ_2+(\lambda y^2+\mu xy^2)dx]=J_1+J_2.$$
Below, all calculations concerning one-forms are performed modulo the forms $dF+GdH$ since they give zero result upon integrating.
We begin with $J_1$. It is convenient to use splitting  $q_{30}= q_{30}^1-m\bar{Q}_1+q_{30}^2$ where the latter collects the integral
terms of $q_{30}$. As $J_1=-\int_{\delta(h)}Q_1dq_{30}$ one needs first to simplify the one-form $dq_{30}$. Thus,
$$\begin{array}{rl}dq_{30}^1 &\sim [\frac{10}{9}q_1^2-(\frac83a_{01}+2a_{03}H+\frac73a_{03}y^2)q_1
+(a_{01}^2+2a_{01}a_{03}y^2+a_{03}^2y^4)]dq_1\\[2mm]
& +[-\frac73a_{03}yq_1^2+(4a_{01}a_{03}y+4a_{03}^2y^3)q_1]dy.
\end{array}$$
Similarly, since $(x^3-x)dx\sim -ydy$ (because $dH\sim 0$) and $(\frac14x^4-\frac12x^2)(x^3-x)dx\sim -\frac12(2yH-y^3)dy$,
one obtains
$$\textstyle dq_{30}^2\sim \varphi dq_1-(6a_{01}a_{03}-2a_{03}b_{10}+12a_{03}^2H)q_1ydy+3a_{03}q_1^2ydy
+6a_{03}^2q_1(2yH-y^3)dy.$$
Therefore,
\begin{equation}\label{E13}
\begin{array}{rl}dq_{30} & \sim[\frac{10}{9}q_1^2-(\frac83a_{01}+2a_{03}H+\frac73a_{03}y^2)q_1
+a_{01}^2+2a_{01}a_{03}y^2+a_{03}^2y^4+\varphi]dq_1\\[2mm]
& -m\bar{Q}'_1dx+[\frac23a_{03}yq_1^2+(2a_{03}b_{10}-2a_{01}a_{03})yq_1-2a_{03}^2y^3q_1]dy.
\end{array}
\end{equation}
Let us turn now to $J_2$. Since $dQ_2=\partial_y Q_2dy+\frac16y^2dq_1$ and $(\lambda+\mu x)y^2d x=-\frac{1}{2}y^2d q_1$,
one can rewrite the initial integral as
$J_2=\int_{\delta(h)}q_{31}[\partial_y Q_2dy-\frac13y^2dq_1]$. To find the impact on $J_2$ of the last term  $2a_{03}q_1Q_1$
in $q_{31}$ (see (\ref{E11})), we note that the first term in $2a_{03}q_1[\partial_y Q_2dy-\frac13y^2dq_1]$ is just the last expression in
(\ref{E13}), while the second term which equals $-\frac23a_{03}y^2q_1dq_1$ will reduce to $\frac53$ the coefficient $\frac73$ in (\ref{E13})
(since $J_1=-\int_{\delta(h)} Q_1dq_{30}$).
So, we will use in what follows $q_{31}$ with its last term removed, and as a compensation,  $q_{30}$  with coefficient $\frac53$
instead of $\frac73$ and its last expression which contains $dy$ removed as well.

With this convention, we obtain by simple direct calculations the expression
$$\textstyle q_{31}[\partial_y Q_2dy-\frac13y^2dq_1]\sim m[-\frac{20}{27}q_1y^3+\frac89a_{01}y^3+\frac23a_{03}y^5]dq_1.$$
On the other hand, using the second equation in (\ref{E10}), we calculate
$$-Q_1dq_1=(2\lambda+2\mu x)Q_1dx=my(x^3-x)dx\sim-my^2dy.$$
Therefore,
$$\begin{array}{l}-Q_1dq_{30}\sim -m[\frac{10}{9}q_1^2-(\frac83a_{01}+2a_{03}H+\frac53a_{03}y^2)q_1
+a_{01}^2+2a_{01}a_{03}y^2+a_{03}^2y^4+\varphi]y^2dy\\[2mm]
 +m\bar{Q}'_1Q_1dx\sim-m[\frac{10}{9}q_1^2y^2-\frac83a_{01}q_1y^2-2a_{03}q_1Hy^2-\frac53a_{03}q_1y^4
+\varphi y^2]dy+m\bar{Q}'_1Q_1dx\\[2mm]
\sim m[\frac{20}{27}q_1y^3-\frac89a_{01}y^3-\frac23a_{03}Hy^3-\frac13a_{03}y^5]dq_1+\frac13my^3d\varphi+m\bar{Q}'_1Q_1dx.
\end{array}$$
Taking into account all terms to integrate coming from $J_1$ and $J_2$, we see that the one-form should be
$q_3\omega\sim m(\frac13a_{03}y^5-\frac23a_{03}Hy^3)dq_1+\frac13my^3d\varphi+m\bar{Q}'_1Q_1dx.$
By the formula of $d\varphi$ from the beginning of the proof,  everything that remained still to integrate is
$M_4(h)=\int_{\delta(h)}m\bar{Q}'_1Q_1dx$ which yields the result stated in (i).

\vspace{1ex}
(ii)  Under the conditions of Proposition 4 (i), system (\ref{E2}) which satisfies equations $M_1=M_2=M_3=M_4=0$ is reduced to
$$\begin{array}{l}
\dot{x}=H_y+\varepsilon[a_{00}(1-x^2)+a_{01}y +a_{21}x^2y+a_{03}y^3],\\[2mm]
\dot{y}=-H_x+\varepsilon[b_{10}(x-x^3)+2a_{00}xy+2a_{21}xy^2].
\end{array}$$

If $a_{21}=0$, then the system is Hamiltonian. Otherwise, the general first integral of the system is of  Darboux type,

$$H_\varepsilon=(x^2+A_1y^2+B_1y+C_1)^{n_1}(x^2+A_2y^2+B_2y+C_2)^{n_2},$$

where
$$\begin{array}{l}{\displaystyle A_{1,2}=\frac{\varepsilon(-a_{21}\pm\sqrt\Delta)}{2(1+\varepsilon b_{10})},\quad
 B_1=B_2=-\frac{2\varepsilon a_{00}}{1+\varepsilon b_{10}}, \quad
 n_{1,2}=\left(1\pm\frac{3a_{21}}{\sqrt\Delta}\right)\frac{1+\varepsilon b_{10}}{4},}\end{array}$$
$$C_1=\frac{1+\varepsilon a_{01}+(n_1+n_2)(2A_1-B_1^2)}{2n_2(A_2-A_1)}, \;
C_2=\frac{1+\varepsilon a_{01}+(n_1+n_2)(2A_2-B_2^2)}{2n_1(A_1-A_2)},$$
and we have denoted $\Delta=a_{21}^2-4a_{03}(b_{10}+1/\varepsilon).$
Assuming $\varepsilon$ to be small and positive, then for $a_{03}>0$ one has $\sqrt\Delta=i\sqrt{-\Delta}$
and therefore $A_{1,2}$, $C_{1,2}$ and $n_{1,2}$ are complex-conjugated. Hence, $H_\varepsilon$ is real-valued.
If $a_{03}<0$, the coefficients and powers are real-valued. However, $C_1<0<C_2$ which implies that
$x^2+A_1y^2+B_1y+C_1< 0$ inside the eight loop. For this reason, we can take in this case
$$H_\varepsilon=(-x^2-A_1y^2-B_1y-C_1)^{-n_1}(x^2+A_2y^2+B_2y+C_2)^{-n_2},$$
as a first integral instead of the former one. Clearly, in both cases $H_\varepsilon$ is analytic first integral
in any bounded and open domain $K\supset \{H\leq 0\}$ as long as $\varepsilon$ is sufficiently small.
Therefore the system has analytic first integral at least in $K$ and the critical points $(\pm 1,0)$ should be centers.

When $a_{03}=0$, the first integral can be obtained by replacing this value in the general formula.
For example, if $a_{21}>0$ one obtains
$$H_\varepsilon=(x^2+B_1y+C_1)^{n_1}(x^2+A_2y^2 +B_2y+C_2)^{n_2},$$
where $n_1=1+\varepsilon b_{10}$,  $n_2=-\frac12(1+\varepsilon b_{10})$,
$A_2=-\frac{\varepsilon a_{21}}{1+\varepsilon b_{10}}$,
$C_1=\frac{1+\varepsilon a_{01}}{\varepsilon a_{21}}-\frac{2\varepsilon a_{00}^2}{a_{21}(1+\varepsilon b_{10})}$,
$C_2=\frac12(C_1-1)$ and $B_{1,2}$ are as above.

Similarly, if $a_{21}<0$, we have
$$H_\varepsilon=(-x^2-B_1y-C_1)^{n_1}(-x^2-A_2y^2 -B_2y-C_2)^{n_2}$$
with the same values as in the case $a_{21}>0$. $\Box$

\vspace{2ex}
\noindent
{\bf Example.} Take for simplicity $\varepsilon=1$ and $a_{03}=0$. Consider the system
$$\begin{array}{l}
\dot{x}=H_y+\varepsilon[1-x^2+10y +2x^2y],\\[2mm]
\dot{y}=-H_x+\varepsilon[x-x^3+2xy+4xy^2]
\end{array}$$
which is neither Hamiltonian nor symmetric with respect to $y$. According to the calculations above its first integral is
$$H_1=\frac{(x^2-y+5)^2}{x^2-y^2-y+2}.$$
The equation of the double loop through the saddle $(0,-\frac{1}{11})$ is
$$(x^2-y+5)^2=\frac{112}{9}(x^2-y^2-y+2).$$
The loop lies in the domain $y<x^2+5$, $x^2-y^2-y+2>0$  together with the centers $(\pm 1,0)$ inside the loop.
The example shows the evolution of the loop while $\varepsilon$ has run along the interval $[0,1]$.


\section  {The zeros of $M_2(h)$ and $M_2^*(h)$  in $\Sigma=(-\frac14,0)$}

In this section we study mainly $M_2(h)$ which could produce most limit cycles among $M_k(h)$. We first recall some known facts about
the Picard-Fuchs system satisfied by the basic integrals $I_k(h)$, $k=0,1,2$ for $h\in\Sigma$ and derive several consequences.
Let us recall that $I_k(h)$ and $I_k'(h)$ are positive in  $\Sigma$ for any $k$ which follows directly from their integral representation
over $\delta(h)$.

\vspace{2ex}
\noindent
{\bf Lemma 1.} (i) {\it The integral $I_1(h)$ satisfies equation $I_1=(h+\frac14)I_1'$.
Hence  $I_1(h)=4I_1(0)(h+\frac14)=\pi\sqrt{2}(h+\frac14)$.}

\vspace{2ex}
\noindent
(ii) {\it The couple ${\textbf{I}}(h)=(I_0(h), I_2(h))^\top$ satisfies system
\begin{equation}\label{PF}
{\textbf{I}}(h)={\textbf{A}}(h){ \textbf{I}}'(h), \qquad {\textbf{A}}(h)=\left(\begin{array}{cc} \frac43h & \frac13\\[2mm]
\frac{4}{15}h& \frac{4}{15}(3h+1)\end{array}\right), \;\;h\in\Sigma.
\end{equation}}

\vspace{2ex}
\noindent
{\bf Proof.} We note that by (\ref{E1}), on each oval $\delta(h)$ the following identities hold
\begin{equation}\label{E15}
y^2=2h+x^2-\frac{x^4}{2},\qquad ydy=(x-x^3)dx, \qquad  y\frac{dy}{dh}=1.
\end{equation}
Therefore for any $k$
$$I_k=\int_{\delta(h)} x^k ydx=\int_{\delta(h)} \frac{x^k y^2}{y}dx=\int_{\delta(h)}
 \frac{x^k(2h+x^2-\frac{x^4}{2})}{y}dx=2h I_k'+I_{k+2}'-\frac{1}{2}I_{k+4}'. $$
Also, after integrating by parts, the second identity in (\ref{E15}) implies
\begin{equation}\label{E16}
I_k=-\frac{1}{k+1}\int_{\delta(h)} x^{k+1}dy=\frac{1}{k+1}\int_{\delta(h)} \frac{x^{k+1}(x^3-x)}{y}dx
=\frac{I_{k+4}'-I_{k+2}'}{k+1}.
\end{equation}
We eliminate $I_{k+4}'$ to obtain $I_k=\frac{1}{k+3}(4hI_k'+I_{k+2}')$. Taking $k=0,1,2$, one derives the system
\begin{equation}\label{E17}
\begin{array}{l}
I_0=\frac{4}{3}hI_0'+\frac{1}{3}I_{2}',\\[2mm]
I_1=hI_1'+\frac{1}{4}I_{3}',\\[2mm]
I_2=\frac{4}{5}hI_2'+\frac{1}{5}I_{4}'.
\end{array}
\end{equation}
The second equation in (\ref{E15}) implies that $I_1(h)= I_3(h)$ which proves (i). With $k=0$,  (\ref{E16}) becomes $I_4'=I_0+I_2'$.
Replacing in the last equation of the system and making use of its first equation too, we obtain the needed form of last equation in
(\ref{E17}) :
$$I_2=\frac{4}{15}hI_0'+\frac{4}{15}(3h+1)I_{2}'$$
which proves statement (ii) in Lemma 1. $\Box$

\vspace{2ex}
\noindent
{\bf Corollary 1.} {\it  With $s=h+\frac14$, the following expansions hold with some positive $c_1$:}
$$\begin{array}{l}
I_0=c_1\left(s+\frac38 s^2+\frac{35}{64}s^3+\frac{1155}{1024}s^4+\frac{45045}{16384}s^5
+\frac{969969}{131072}s^6+...\right)\\[2mm]
I_2=c_1\left(s-\frac18 s^2-\frac{5}{64}s^3-\frac{105}{1024}s^4-\frac{3003}{16384}s^5
 -\frac{51051}{131072}s^6-...\right)
\end{array}
$$
{\it In particular, $I''_0$ and all its derivatives are positive whilst $I_2''$ and all its derivatives are negative in $\Sigma$.}

\vspace{2ex}
\noindent
{\bf Proof.} As $I_0$ and $I_2$ are analytic and vanish at $h=-\frac14$, one can take $I_0=\sum_{k=1}^\infty a_ks^k$,
$I_2=\sum_{k=1}^\infty c_ks^k$ and replace in system (\ref{PF}) with
$$ {\textbf {\it A}}=\left(\begin{array}{cc} \frac{4s-1}{3} & \frac13\\[2mm]\frac{4s-1}{15} & \frac{12s+1}{15}\end{array}\right).$$
Thus we obtain the following simple recursive formulas
$$a_k=(5-4k)c_k, \quad c_{k+1}=\frac{(4k-5)(4k-3)}{4k(k+1)}c_k, \quad k=1,2,...,$$
which imply
$$ a_1=c_1,\quad a_k=\frac{(4k-5)!!}{4^{k-1}(k-1)!k!}c_1,\quad c_k=-\frac{(4k-7)!!}{4^{k-1}(k-1)!k!}c_1,\quad k=2,3,....$$
By the way, since the oval $\delta(h)$ shrinks to the center $C$ as $h\to-\frac14$, this yields that $I_k/I_0\to 1$ when $s\to 0$.
Hence, $c_1=\pi\sqrt{2}$ and $I_1=c_1s$. $\Box$

\vspace{2ex}
By differentiating the Picard-Fuchs system (\ref{PF}), one obtains

\vspace{2ex}
\noindent
{\bf Corollary 2.} {\it The following equalities hold:
$ {\textbf I}^{(k)}(h)={\textbf A}_k(h){\textbf I}^{(k-1)}(h)$, $k=1,2,3,...$, where ${\textbf A}_1={\textbf A}^{-1}$,
${\textbf A}_{k+1}={\textbf A}'_k{\textbf A}_k^{-1}+ {\textbf A}_k$,  $k=1,2,3,.....$
 Explicitly,
$${\textbf A}_k(h)=\frac{{\textbf B}_k(h)}{h(4h+1)},\qquad
 {\textbf B}_1=\left(\begin{array}{cc}3h+1 & -\frac54\\[1mm]-h & 5h\end{array}\right), $$
$${\textbf B}_2=\left(\begin{array}{cc} -h & -\frac14\\[1mm]-h & h\end{array}\right),
{\textbf B}_3=\left(\begin{array}{cc}-(5h+1) &\frac34\\[1mm]-h & -3h\end{array}\right),
{\textbf B}_4=\left(\begin{array}{cc}-(9h+2) & \frac74\\[1mm]-h & -7h\end{array}\right).$$}

\vspace{2ex}
Next, we let  $\displaystyle \omega(h)=\frac{I_2(h)}{I_0(h)}$ and  $\displaystyle \nu(h)=\frac{I'''_2(h)}{I'''_0(h)}$, $h\in \Sigma.$
Then the associated to system {\it ${\textbf I}^{(4)}(h)={\textbf A}_4(h){\textbf I}'''(h)$}  Riccati equation for $\nu(h)$ reads
\begin{equation}\label{E18}
\nu'(h)=\frac{\frac74\nu^2-2(h+1)\nu +h}{-h(4h+1)},
\end{equation}
which is equivalent to the system
\begin{equation}\label{E19}
\begin{array}{l}
\frac{d\nu}{dt}=\frac74\nu^2-2(h+1)\nu+h,\\[2mm]
\frac{dh}{dt}=-h(4h+1).
\end{array}
\end{equation}
The phase portrait of the system is shown in Figure 1.

\vspace{2ex}
\noindent
{\bf  Lemma 2.} {\it  The graph of $\nu=\nu(h)$  is the unstable manifold of system $(\ref{E19})$ at the saddle
point $S(-\frac14,-\frac17)$, connecting it to the stable node at the origin. Moreover, $\nu'(h)>0$ and $\nu''(h)>0$
in $\Sigma$.}

\vspace{2ex}
\noindent
{\bf Proof.} Applying Corollary 2 three times, we obtain
{\it $${\textbf I}'''(h)=\frac{{\textbf B}_3(h){\textbf B}_2(h){\textbf B}_1(h)}{h^3(4h+1)^3}{\textbf I}(h)=\frac{15}{16h^2(4h+1)^2}
\left(\begin{array}{cc} 4h & 1\\[1mm] 4h & -4h\end{array}\right){\textbf I}(h).$$}
By this identity, we have
\begin{equation}\label{E20}
\nu(h)=\frac{4h(1-\omega(h))}{4h+\omega(h)}.
\end{equation}
Direct calculation of the respective integrals over $\delta(0)$ yields values $I_0(0)=\frac43$, $I_2(0)=\frac{16}{15}$.
Hence $\omega(0)=\frac45$ which by (\ref{E20}) means that $\nu(h)\to 0$ as $h\to 0$. Next, we can use Corollary 1 to calculate
$\omega=1-\frac12s-\frac{7}{16}s^2-\frac{203}{256}s^3 + ...$ and then use (\ref{E20}) to calculate
$\nu=-\frac17+\frac37s+\frac{33}{112}s^2+...$ in order to verify that $\nu(-\frac14)=-\frac17$, $\nu'(-\frac14)=\frac37$
and $\nu''(-\frac14)>0$. On the other hand, the horizontal isocline $\nu=\psi_1(h)$  is the lower branch of the hyperbola
$\frac74\nu^2-2(h+1)\nu+h=0$ in the $(h,\nu)$-plane, going also through $S$ and $O$, see Figure 1A. Differentiating the equation with respect
to $h$, we obtain that $\psi_1'(-\frac14)=\frac{9}{14}>\nu'(-\frac14)$. This implies that the curve $\nu=\nu(h)$ always stays below the
horizontal isocline, hence $\nu'(h)>0$.

\begin{figure}[!htbp]
\begin{center}
\psfrag{1}{$\nu=\nu(h)$}\psfrag{2}{$\nu=\psi_1(h)$}\psfrag{3}{$\nu=\psi_2(h)$}\psfrag{H}{$h$}
\psfrag{V}{$\nu$}\psfrag{O}{{\small$O$}}
\psfrag{A}{$-\frac{1}{4}$}\psfrag{S}{$S$}\psfrag{4}{(A)}\psfrag{5}{(B)}\psfrag{6}{$h_0$}\psfrag{7}{$h_1$}
\psfrag{8}{$h_2$}\psfrag{9}{$l$}
\psfig{file=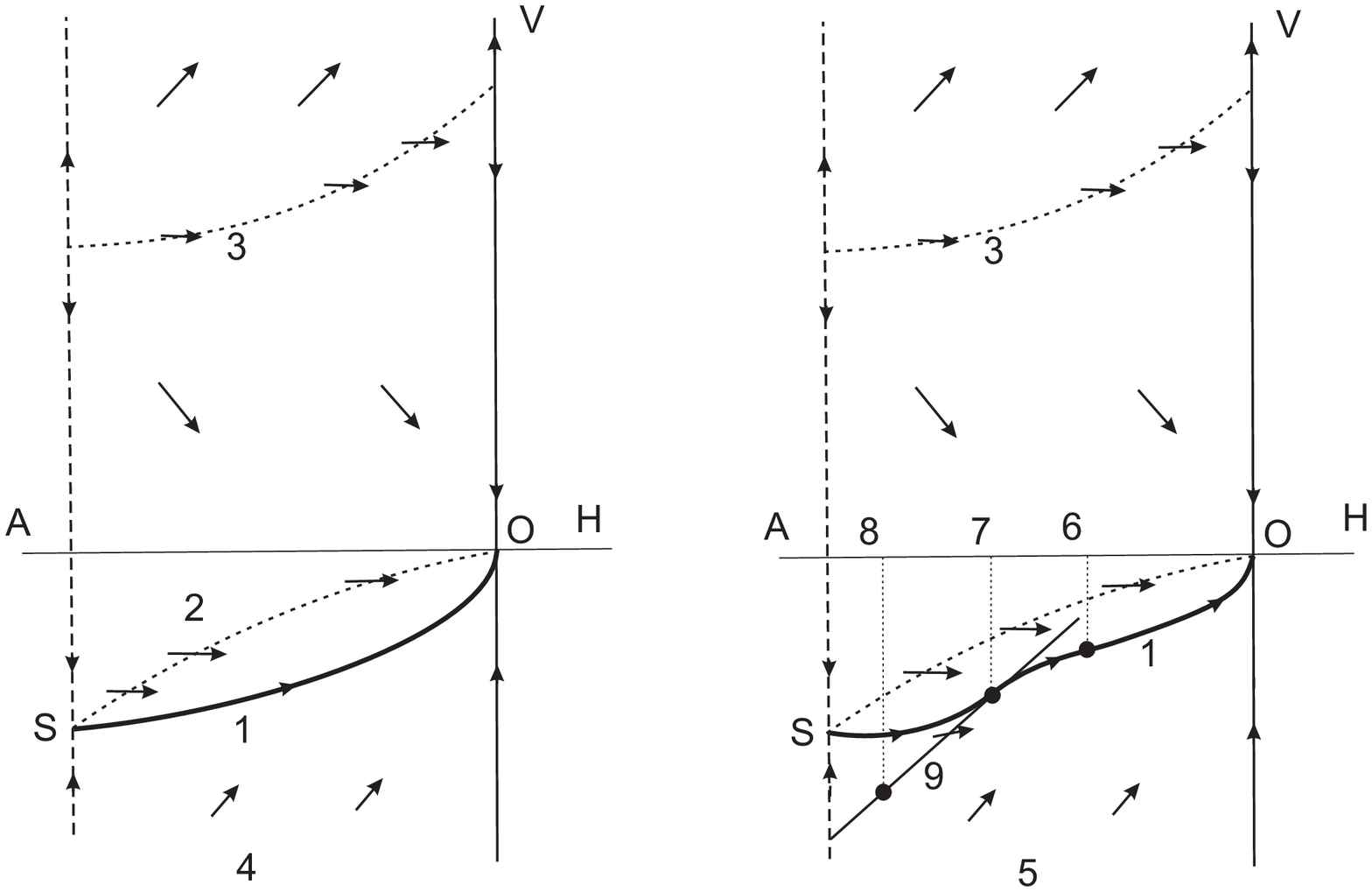,height=75mm,width=120mm}
\centerline{\footnotesize Figure 1. The phase portraits of system (\ref{E19})}
\end{center}
\end{figure}

It remains to verify inequality $\nu''(h)>0$.
Assume that there is $h_0\in \Sigma$ such that $\nu''(h_0)\leq 0$. Let $h_1\leq h_0$ be the unique value such that
$\nu''(h_1)=0$ and $\nu''(h)>0$ in $(-\frac14, h_1)$.  Consider the tangent line $\ell: \nu=L(h)=\nu'(h_1)(h-h_1)+\nu(h_1)$
having at $(h_1,\nu(h_1))$ a tangent point of multiplicity at least two. By $\nu'(h_1)>0$ and convexity of $\nu(h)$ in
$(-\frac14, h_1)$, $\ell$ intersects the separatrix $\{h=-\frac14\}$ below $S$, see Figure 1B.
By the saddle property, there is on $\ell$ another contact point with the field at $h_2\in (-\frac14,h_1)$.
This is however impossible  because  by using (\ref{E19}), we see that the expression
$$\left.\frac{d\nu}{dt}-L'(h)\frac{dh}{dt}\right|_{\nu=L(h)}$$
is a polynomial in $h$ of degree 2. Hence, $\nu''>0$ in $\Sigma$. $\Box$

\vspace{2ex}
Now we are prepared to establish the following proposition which proves Theorem 2:

\vspace{2ex}
\noindent
{\bf  Proposition 5.} {\it  The second Melnikov functions $M_2(h)$ and $M^*_2(h)$ can have each at most 5 zeros in $\Sigma = (-\frac14,0)$,
counted with their multiplicities.}

\vspace{2ex}
\noindent
{\bf  Proof.} We begin by estimating the number of zeros of $M_2'''(h)$. We wish to prove that $M_2'''(h)$ has no more that 3 zeros in $\Sigma$.
By Proposition 2(i) and Lemma 1(i), we have
$$M_2'''(h)=(\alpha_0+\alpha_1 h)I_0'''(h)+(\gamma_0+\gamma_1 h)I_2'''(h)+3 \alpha_1 I_0''(h)+3 \gamma_1 I_2''(h).$$
 From Corollary 2 one obtains
{\it $${\textbf  I}''(h)={\textbf A}_3^{-1}(h){\textbf  I}'''(h)=\left(\begin{array}{cc} -\frac45h & -\frac15\\[2mm]
 \frac{4}{15}h  &  -\frac43h - \frac{4}{15} \end{array}\right){\textbf  I}'''(h).$$}
Therefore
$$M_2'''(h)=(a h+b)I_0'''(h)-(ch+d) I_2'''(h)$$
where $a,b,c,d$ are independently free constants calculated explicitly from $\alpha_i$ and $\gamma_i$.
Obviously, if all the constants are zero, then $M_2$ has no more than two zeros, because by Theorem 1 we suppose that $M_2$
is not identically zero. By Lemma 2, if $ad-bc=0$ then $M_2'''$
has at most two zeros in $\Sigma$; the same holds if $c=0$. Below we consider the remaining cases when $c(ad-bc)\neq 0$.
If $h=-d/c$, then $M_2'''(h)\neq 0$ so we can without any loss of zeros in $\Sigma$ to suppose $h\neq -d/c$ and to rewrite $M_2'''(h)$ as
$$M_2'''(h)=(ch+d)I_0'''(h)(N(h)-\nu(h)),\quad N(h)=\frac{ah+b}{ch+d}.$$
We  have to determine the number of intersections between the hyperbola $\nu=N(h)$ and the separatrix trajectory $\nu=\nu(h)$.
If $ad-bc<0$ then $N'<0$ and there are at most two intersections, see Figure 2B. If $ad-bc>0$, then the separatrix could intersect only
one of the branches. At most two intersections are possible with the right (lower) branch because of its concavity.

\begin{figure}[!htbp]
\begin{center}
\psfrag{A}{\small {Case (A): $ad-bc>0$}}\psfrag{B}{\small {Case (B): $ad-bc<0$}}
\psfig{file=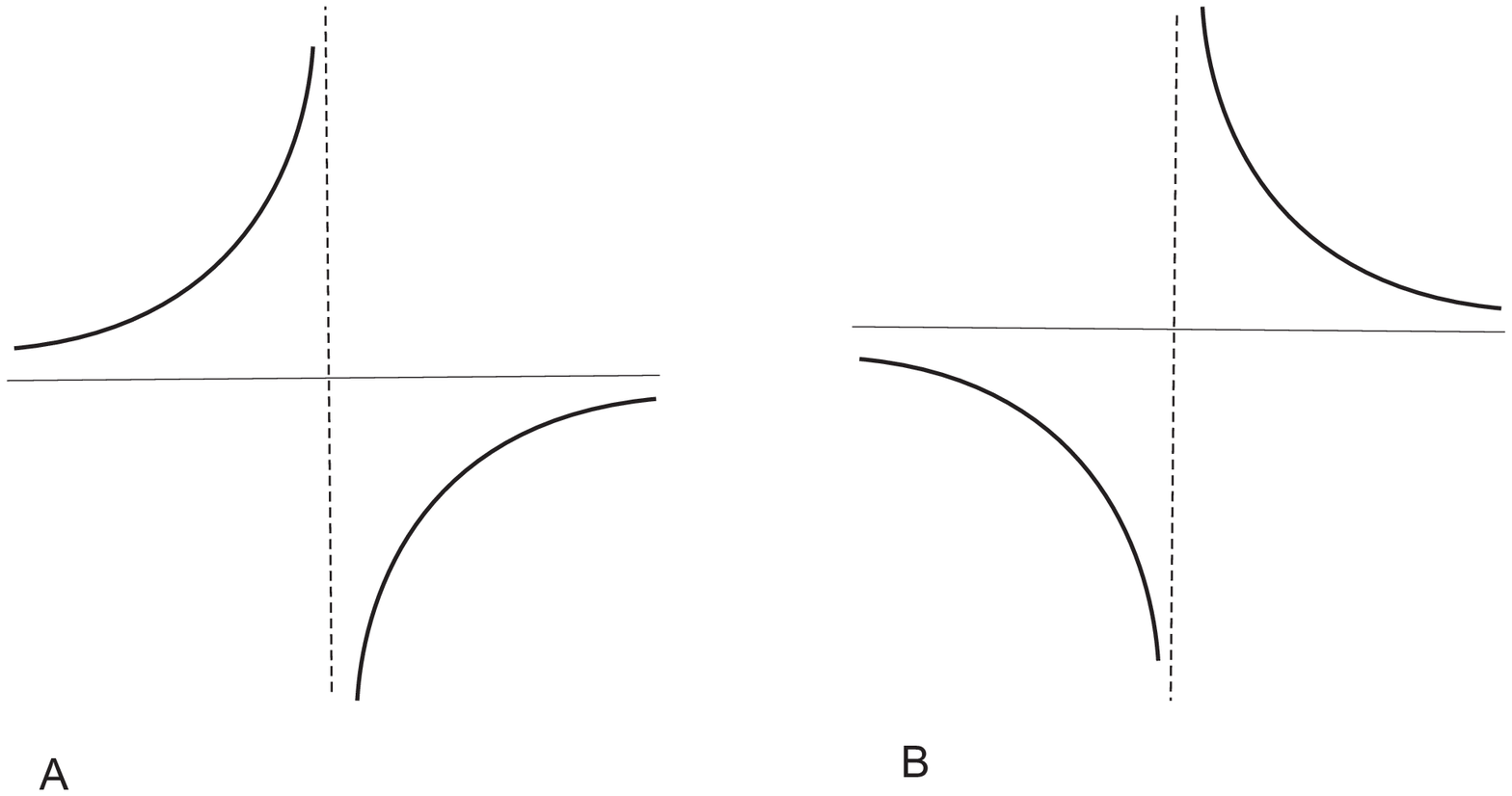,height=50mm,width=90mm} \centerline{\footnotesize Figure 2. The behavior of the curve $\nu=N(h)$.}
\end{center}
\end{figure}

In the left branch case, by using the saddle property at $S$, on the left of the first intersection point between these two curves there is a contact
point on $\nu=N(h)$ with the vector field  (\ref{E19}). Hence the number of intersection points of these two curves is controlled by the
number of tangent points on $\nu=N(h)$ with respect to the vector field (\ref{E19}). This is true even if the hyperbola goes through
the critical point $S$ or $O$ or both, since then $S$ (respectively $O$) is simultaneously an intersection and contact point.
Noticing that the numerator of
$$\left.\frac{d\nu}{dt}-N'(h)\frac{dh}{dt}\right|_{\nu=N(h)}$$
is a polynomial in $h$ of degree 3, we deduce that the number of zeros of $M_2'''(h)$ is also at most 3, including their multiplicities.

Therefore the number of zeros of $M_2(h)$ is at most 6, including their multiplicities. Since $M_2(-\frac14)=0$, this number
 in $h\in (-\frac14,0)$ is at most 5. As our proof above uses information about $M_2'''=M_2^{*'''}$ only, the same conclusion holds for
$M_2^*(h)$, too. $\Box$


\section  {The distribution of limit cycles on both nests}

\vspace{2ex}
In this section we investigate the simultaneous bifurcation of limit cycles from both period annuli $\mathcal{A}$, $\mathcal{A}^*$ in order to prove
claims (i)-(iii) of Theorem 3. First of all, one can proceed as above to obtain the following result, which implies (iii):

\vspace{2ex}
\noindent
{\bf  Proposition 6.} {\it  The particular perturbations with
$$ M_2(h)= M^*_2(h)=(\alpha_0+\alpha_1h)I_0(h)+(\gamma_0+\gamma_1h)I_2(h),\quad h\in \Sigma $$
produce distribution $(m,m)$ of limit cycles where $m\leq 3$.}

\vspace{2ex}
\noindent
{\bf  Proof.} The proof goes just as the proof of  Proposition 5 with $\omega(h)=\frac{I_2(h)}{I_0(h)}$ replacing $\nu(h)$
and another hyperbola in the $(h,\omega)$-plane given by $\omega=N(h)=-\frac{\alpha_1h+\alpha_0}{\gamma_1h+\gamma_0}$.
Thus, one can rewrite the above formula as $M_2(h)=(\gamma_1h+\gamma_0)I_0(h)[\omega(h)-N(h)]$, $h\in \Sigma.$
The Riccati equation for $\omega(h)$ related to system (\ref{PF}) reads
$$\omega'(h)=\frac{\frac54\omega^2+(2h-1)\omega-h}{h(4h+1)}, \quad h\in \Sigma,$$
which is equivalent to the system
\begin{equation}\label{E21}
\begin{array}{l}
\frac{d\omega}{dt}=-\frac54\omega^2-(2h-1)\omega+h,\\[2mm]
\frac{dh}{dt}=-h(4h+1).
\end{array}
\end{equation}
The phase portrait of the system is shown in Figure 3.

\begin{figure}[!htbp]
\begin{center}
\psfrag{1}{$\omega=\omega(h)$}\psfrag{2}{$\omega=\varphi_1(h)$}\psfrag{3}{$\omega=\varphi_2(h)$}
\psfrag{H}{$h$}\psfrag{W}{$\omega$}\psfrag{O}{{\small$O$}}
\psfrag{A}{$-\frac{1}{4}$}\psfrag{S}{$S$}\psfrag{B}{$N$}\psfrag{C}{$M$}
\psfig{file=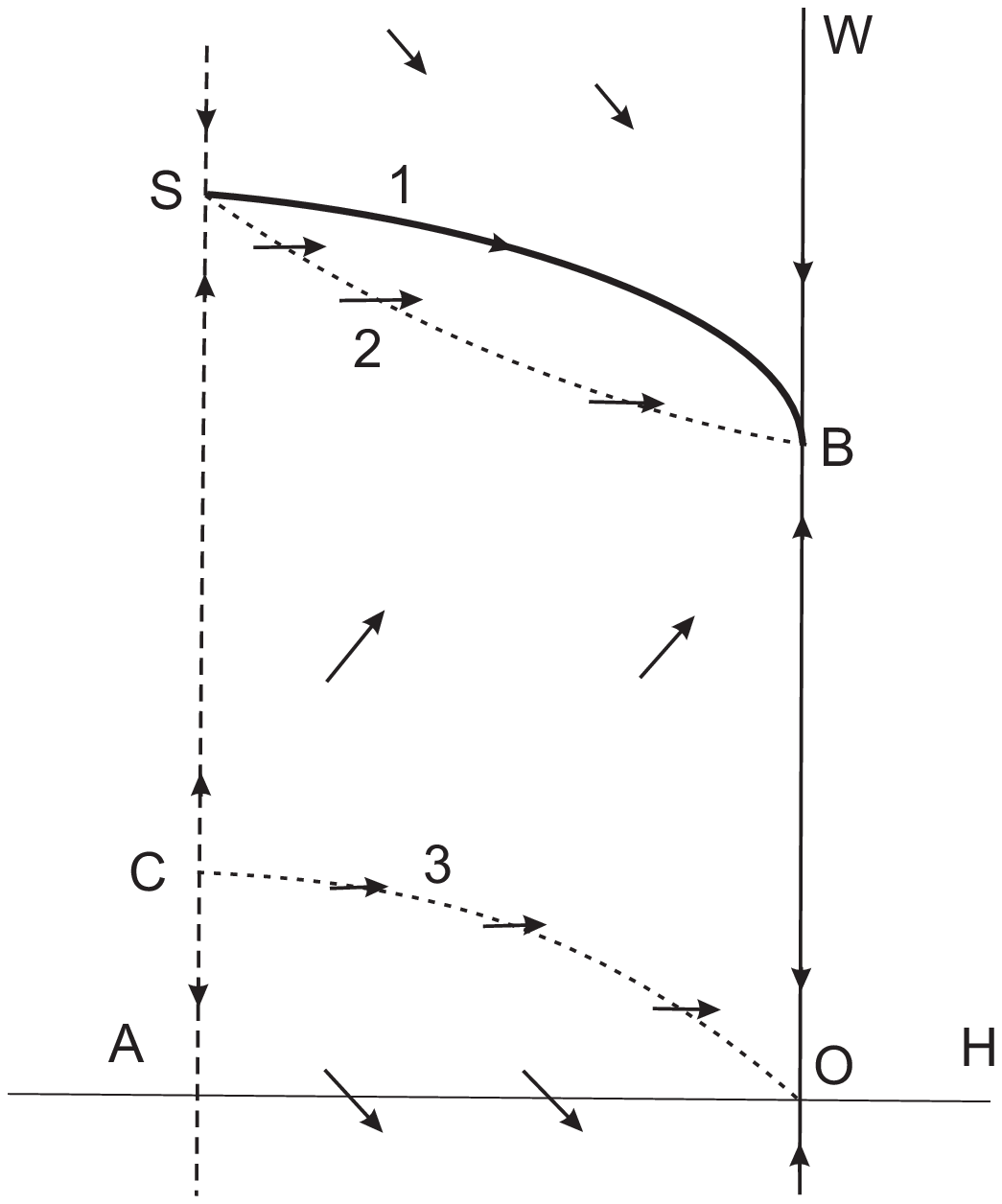,height=75mm,width=80mm} \centerline{\footnotesize
Figure 3. The phase portrait of system (\ref{E21})}
\end{center}
\end{figure}

We already have established that the graph of $\omega(h)$ coincides with the separatrix trajectory connecting the saddle $S(-\frac14,1)$
and the stable node $N(0,\frac45)$. Moreover, $\omega'(-\frac14)=-\frac12$ and $\omega''(-\frac14)<0$. On the other hand, the upper
branch $\omega=\varphi_1(h)$ of the hyperbola $-\frac54\omega^2-(2h-1)\omega+h=0$ which is the horizontal isocline going through
$S$ and $N$ satisfies $\varphi'_1(-\frac14)=-1$. Then exactly as in Lemma 2, one can prove that $\omega'(h)<0$,  $\omega''(h)<0$,
$h\in\Sigma$, namely $\omega=\omega(h)$ is strictly decreasing and strictly concave. As a result, there are at most 3 intersections
with the hyperbola, the most delicate case occurs in handling the lower (left) branch in Figure 2B.  $\Box$

\vspace{2ex}
\noindent
{\bf Proof of Theorem 3 (i),(ii).}  Instead of $M_2(h)$ and $M_2^*(h)$, let us consider their ratios with $I_1(h)$, rewritten as
$$R_\pm(h)=(\alpha_0+\alpha_1h)J_0(h)\pm(\beta_0+\beta_1h)+(\gamma_0+\gamma_1s)J_2(h)=R_0(h)\pm(\beta_0+\beta_1h),$$
where $J_0(h)=I_0/I_1$, $J_2(h)=I_2/I_1$ and $h\in\Sigma=(-\frac14,0)$. Clearly, $M_2(h)$ and $R_+(h)$ have the same number
of zeros in $\Sigma$ and the same holds for $M_2^*(h)$ and $R_-(h)$.

Then the general question about the total number of zeros of both $M_2$, $M_2^*$ in $\Sigma$ could be investigated by counting the
intersection points (taken with multiplicities) of the graph of $R_0(h)$ with both lines $r=\pm(\beta_0+\beta_1h)$ in the $(h, r)$-space,
as shown in Figures 4, 5 and 6, or equivalently, the intersections of the line $r=\beta_0+\beta_1h$ with both mirror curves $r=\pm R_0(h)$,
or finally, the intersection points between the graphs of  $|R_0(h)|$ and  $|\beta_0+\beta_1h|$ in the half-plane $r\geq 0$.

Note that the intersections on $r=0$ should be counted twice.

Next, consider the ratio $R_0(h)$, $h\in \Sigma$. We claim that:

\vspace{1ex}

1). The function $R_0(h)$ has at most 3 zeros in $\Sigma$.

2). The graph of $R_0(h)$ has at most 3 extrema in $\Sigma$.

3). The graph of $R_0(h)$ has at most 3 inflexion points in $\Sigma$.

\vspace{1ex}
 In fact, we claim above that each of the equations $R_0(h)=0$,  $R'_0(h)=0$, $R''_0(h)=0$ has at most 3 zeros in $\Sigma$.
To justify the geometry of the curve $R_0(h)=0$, we will also need its asymptotic behavior near the endpoints of $\Sigma$.
Consider first $h=-0$. By using the first system in (\ref{22}) below, one can easily verify that
the fundamental system of solutions at this end is formed by
$$\left(\begin{array}{c}\varphi_0\\[1mm]\varphi_2\end{array}\right)=\left(\begin{array}{c}-h+\frac{35}{8}h^2-\frac{1155}{64}h^3+...\\[1mm]
\frac12h^2-\frac{21}{8}h^3+\frac{3003}{256}h^4+...\end{array}\right),$$
$$\left(\begin{array}{c}\psi_0\\[1mm]\psi_2\end{array}\right)=
\left(\begin{array}{c} \frac43-\frac{11}{6}h+\frac{179}{24}h^2-\frac{11639}{384}h^3+...  \\[1mm]
\frac{16}{15}-\frac{4}{15}h+\frac{16}{15}h^2-\frac{143}{30}h^3+\frac{313009}{15360}h^4+...
\end{array}\right)+\left(\begin{array}{c}\varphi_0\\[1mm]\varphi_2\end{array}\right)\log(-h).$$
Then
$$\left(\begin{array}{c}J_0\\[1mm]J_2\end{array}\right)=\sigma \left(\begin{array}{c}\psi_0\\[1mm]\psi_2\end{array}\right)
+const\left(\begin{array}{c}\varphi_0\\[1mm]\varphi_2\end{array}\right)$$
where the constant $\sigma = 2\sqrt{2}/\pi$ is determined from direct calculation of the integrals $I_k(0)$, $k=0,1,2.$
This yields immediately that  near $h=-0$, we have
$$\textstyle R_0(0)=\sigma[\frac43\alpha_0+\frac{16}{15}\gamma_0],\;\; R_0'(h)\sim -\sigma\alpha_0\ln(-h),\;\; R_0''(h)\sim
-{\displaystyle\frac{\sigma\alpha_0}{h}}.$$

On the other hand, by the formulas derived in the proof ot Theorem 4 below (with $\beta_0=\beta_1=0$), one obtains
$R_0(-\frac14)=m_1$, $R'_0(-\frac14)=m_2$, $R''_0(-\frac14)=2m_3$,  Hence
$$\begin{array}{l}
 R_0(-\frac14)=\alpha_0+\gamma_0-\frac14\alpha_1-\frac14\gamma_1, \\[2mm]
R'_0(-\frac14)=\frac38\alpha_0-\frac18\gamma_0+\frac{29}{32}\alpha_1-\frac{33}{32}\gamma_1,\\[2mm]
R''_0(-\frac14)=\frac{35}{32}\alpha_0-\frac{5}{32}\gamma_0+\frac{61}{128}\alpha_1-\frac{27}{128}\gamma_1.
\end{array}
$$

Provided 1), 2), 3) are established, then 1) and 2) would imply Theorem 3 (ii) (see Figure 4 below). Similarly, 1), 2) and 3) would
imply Theorem 3 (i) (see Figures 4, 5 and 6 below).

In order to verify 1), 2), 3), let us denote {\it ${\textbf J}=(J_0,J_2)^\top=\left(\frac{I_0}{I_1}, \frac{I_2}{I_1}\right)^\top.$}
It turns out that  ${\textbf {\it J}}$ satisfies useful Picard-Fuchs systems and respective Riccati equations, namely
{\it \begin{equation}\label{22}
\begin{array}{lll}
 {\textbf J}=\left(\begin{array}{cc} -4h & -5\\[1mm] -4h & 4h-4\end{array}\right){\textbf J}', &
{\displaystyle  u'=\frac{\frac54u^2+(2h-1)u-h}{h(4h+1)},} &  {\displaystyle u=\frac{J_2}{J_0},}\\[6mm]
{\textbf J}'=\left(\begin{array}{cc} -\frac45h & -1\\[1mm] \frac{4}{15}h & -\frac43h\end{array}\right){\textbf J}'', &
 {\displaystyle v'=\frac{-\frac{15}{4}v^2+2hv-h}{h(4h+1)},}  & {\displaystyle v=\frac{J_2'}{J_0'},}\\[6mm]
{\textbf J}''=\left(\begin{array}{cc} -\frac49h & -\frac59\\[1mm] \frac{4}{63}h & -\frac47h-\frac{4}{63}\end{array}\right){\textbf J}''', &
 {\displaystyle w'=\frac{-\frac{35}{4}w^2+(2h+1)w-h}{h(4h+1)},} &  {\displaystyle  w=\frac{J_2''}{J_0''}.}
\end{array}
\end{equation}}

\begin{figure}[!htbp]
\begin{center}
\psfrag{R}{$R_0(h)$}\psfrag{S}{$h$}\psfrag{O}{{\small $-\frac{1}{4}$}}\psfrag{1}{$(+)$}
\psfrag{A}{$(-)$}\psfrag{4}{$O$}\psfrag{2}{$(+)$}\psfrag{B}{$(-)$}\psfrag{3}{$(+)$}
\psfrag{C}{$(-)$}\psfrag{E}{$\alpha_0>0, \, 5\alpha_0+4\gamma_0<0$}
\psfig{file=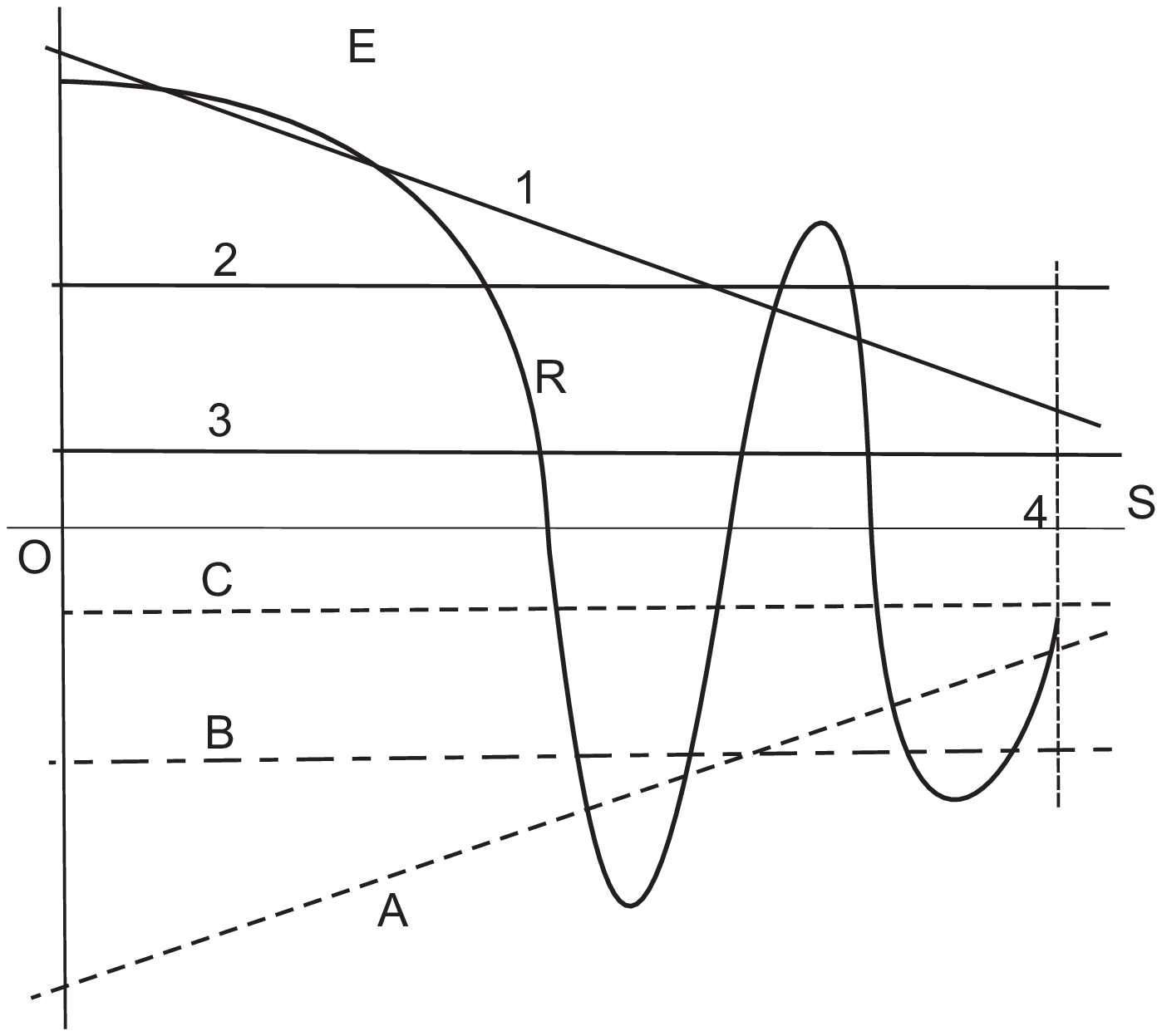,height=60mm,width=75mm}
\centerline{\footnotesize Fig.\,4.  Distributions: $4+4$, $R_0(h)=\pm(\beta_0+\beta_1h)$;  $3+4$, $3+3$,  $R_0(h)=\pm\beta_0$.}
\end{center}
\end{figure}
\begin{figure}[!htbp]
\begin{center}
\psfrag{R}{$R_0(h)$}\psfrag{S}{$h$}\psfrag{O}{{\small $-\frac{1}{4}$}}\psfrag{1}{$(+)$}
\psfrag{A}{$(-)$}\psfrag{4}{$O$}
\psfrag{E}{$\alpha_0>0, \, 5\alpha_0+4\gamma_0<0$}
\psfig{file=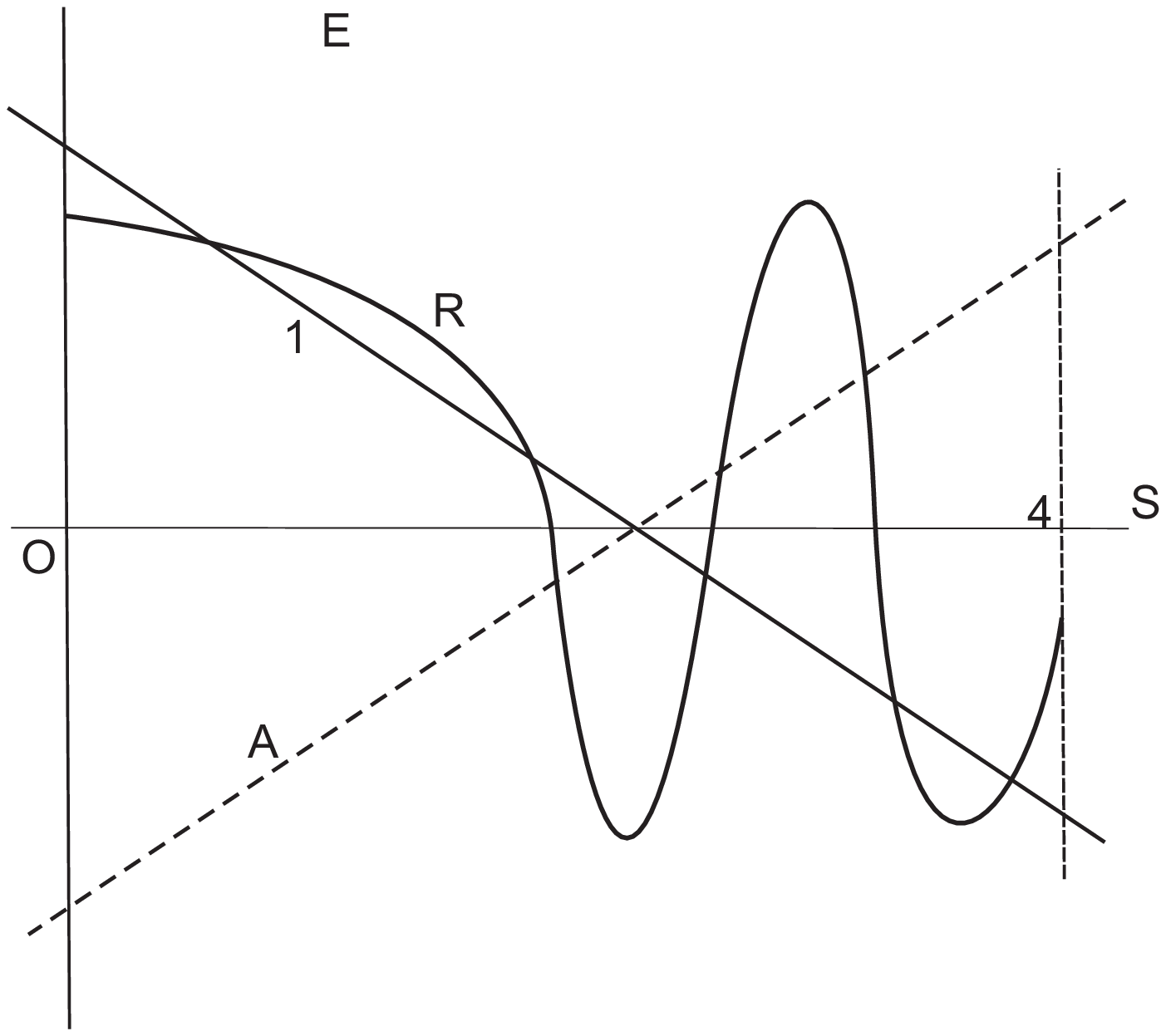,height=60mm,width=75mm}
\centerline{\footnotesize Fig.\,5. Distribution $5+3$,  $R_0(h)=\pm(\beta_0+\beta_1h)$.}
\end{center}
\end{figure}
\begin{figure}[!htbp]
\begin{center}
\psfrag{R}{$R_0(h)$}\psfrag{S}{$h$}\psfrag{O}{{\small $-\frac{1}{4}$}}\psfrag{1}{$(+)$}
\psfrag{A}{$(-)$}\psfrag{4}{$O$}
\psfrag{E}{$\alpha_0>0, \, 5\alpha_0+4\gamma_0>0$}
\psfig{file=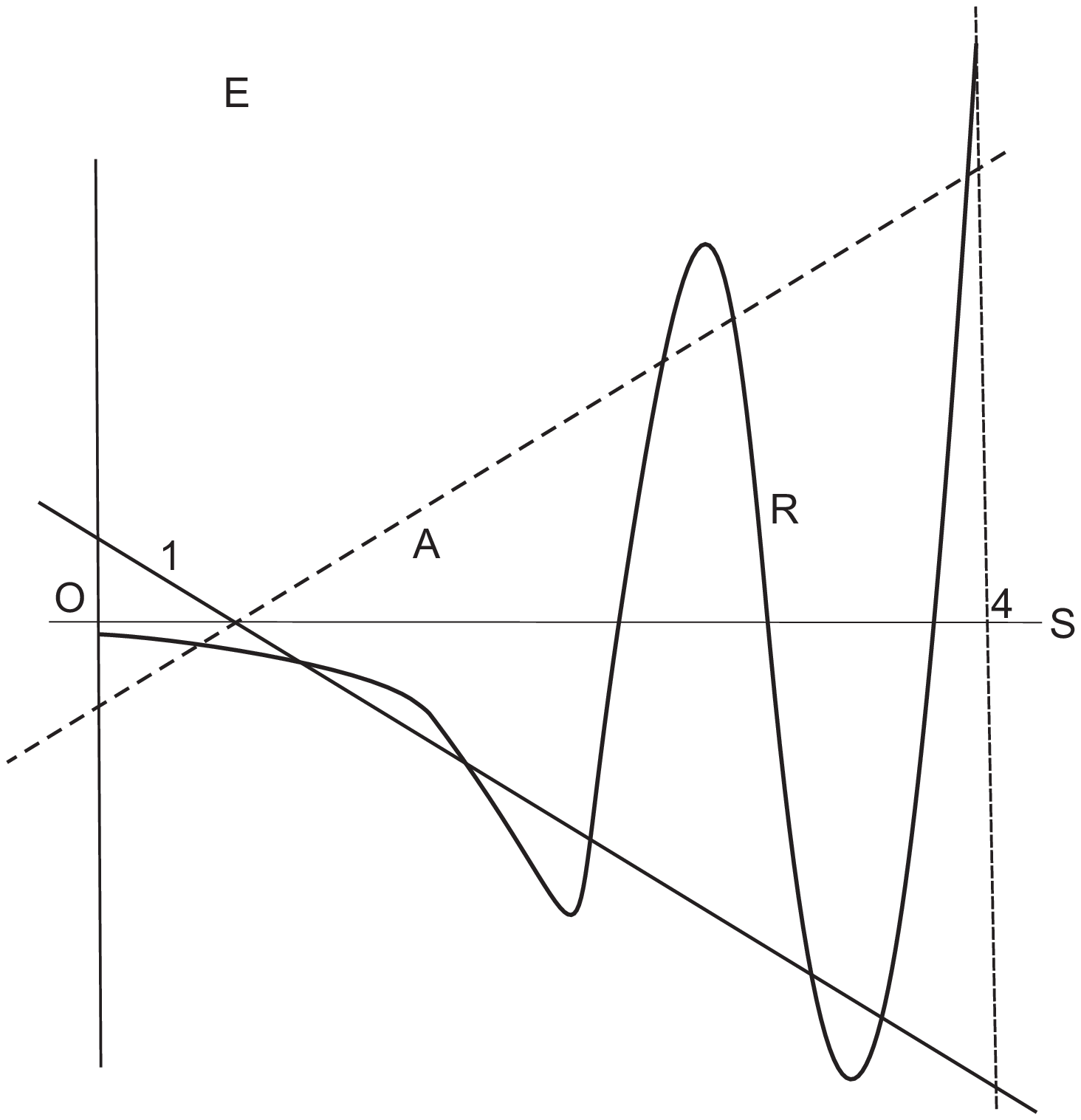,height=60mm,width=75mm}
\centerline{\footnotesize Fig.\,6. Distribution $5+4$, $R_0(h)=\pm(\beta_0+\beta_1h)$.}
\end{center}
\end{figure}

One can derive the first system above by replacing {\it ${\textbf I}={\textbf J}I_1$ in (\ref{PF})}. Then the second and the third ones are obtained
by differentiation. Next, we differentiate $R_0(h)$ and use these systems to obtain the expressions
\begin{equation}\label{23}
R_0^{(k)}(h)=(a_kh+b_k)J_0^{(k)}(h)+(c_kh+d_k)J_2^{(k)}(h),\;\; k=1,2
\end{equation}
with appropriate coefficients $a_k$, etc.

We note that Condition 1) is fulfilled because the ratio $u=J_2/J_0=\omega$ satisfies the same Riccati equation as $\omega$ in the proof of
Proposition 6 above and has the same boundary values as $\omega$.

To establish 2), we use (\ref{23}) with $k=1$ and the system, equivalent to the second Riccati equation in (\ref{22}):
$$
\begin{array}{l}
\frac{dv}{dt}=\frac{15}{4}v^2-2hv+h,\\[2mm]
\frac{dh}{dt}=-h(4h+1).
\end{array}
$$
The phase portrait of the system is shown in Figure 7.

\begin{figure}[!htbp]
\begin{center}
\psfrag{1}{$v=v(h)$}\psfrag{H}{$h$}\psfrag{V}{$v$}\psfrag{O}{{\small$SN$}}\psfrag{A}{{\small $N$}}\psfrag{S}{{\small $S$}}
\psfig{file=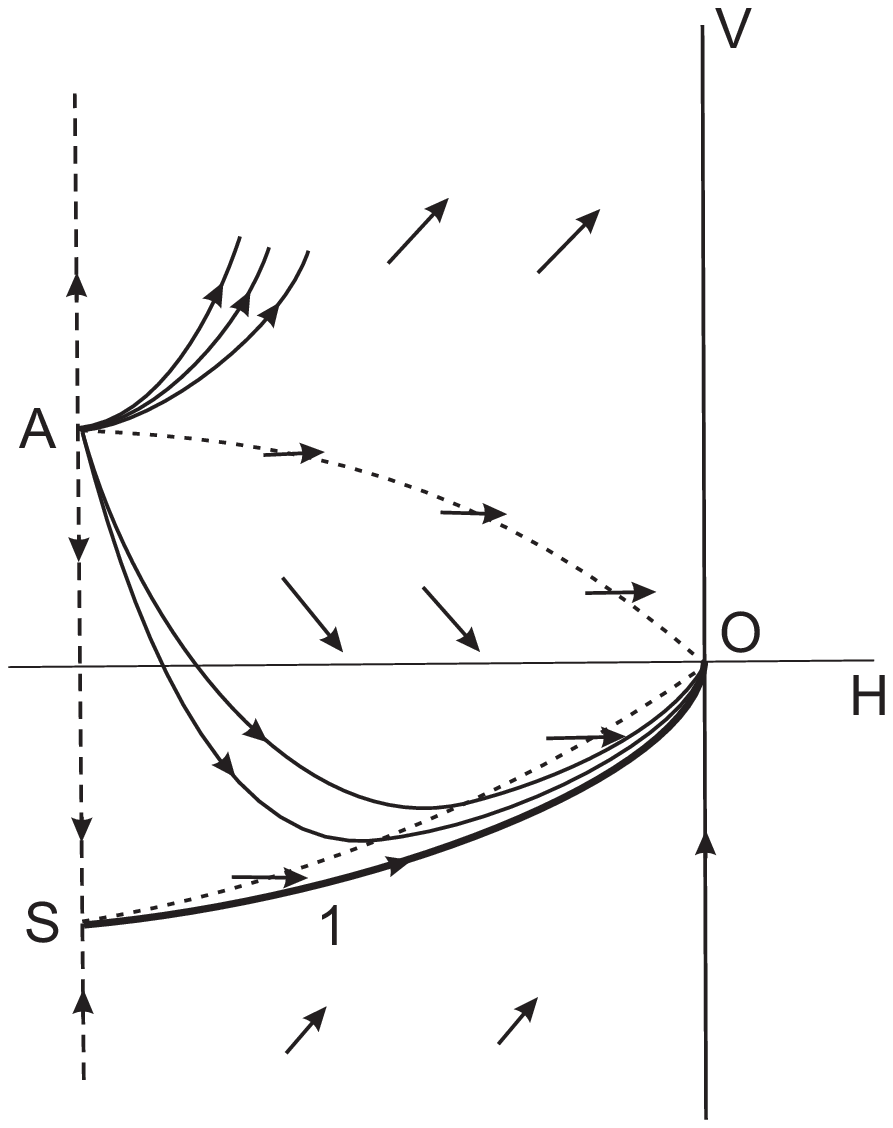,height=80mm,width=60mm} \centerline{\footnotesize Figure \,7. The behavior of $v=v(h)$.}
\end{center}
\end{figure}

One can check that the graph of $v(h)$ coincides with the separatrix connecting the saddle
$S=(-\frac14, -\frac13)$  and the saddle-node $SN$ at the origin. Indeed by calculation using Corollary 1
one obtains  $v=-\frac13+\frac59s+\frac{245}{432}s^2+... $, hence $v(-\frac14)=-\frac13$, $v'(-\frac14)=\frac59$ and
$v''(-\frac14)>0$. The horizontal isocline on the picture is the left branch of the hyperbola $\frac{15}{4}v^2-2hv+h=0$
connecting the three singular points $S$, $SN$ and the node $N(-\frac14,\frac15)$. As its lower part $v=\psi(h)$, connecting
$S$ and $SN$ is increasing and $\psi'(-\frac14)=\frac56>v'(-\frac14)$, we conclude that $v(h)$ is increasing  and convex in
the whole interval and $v(h)<\psi(h)$. The proof is the same as in Lemma 2. Then the proof of 2) is the same as the proof of
Propositions 5 and 6.

Finally, to establish 3), we take (\ref{23}) with $k=2$ and the system, equivalent to the third Riccati equation in (\ref{22}):
$$
\begin{array}{l}
\frac{dw}{dt}=\frac{35}{4}w^2-(2h+1)w+h,\\[2mm]
\frac{dh}{dt}=-h(4h+1).
\end{array}
$$
The phase portrait of the system is shown in Figure 8.

\begin{figure}[!htbp]
\begin{center}
\psfrag{1}{$w=w(h)$}\psfrag{H}{$h$}\psfrag{V}{$w$}\psfrag{O}{{\small$N$}}
\psfrag{A}{{\small ${\bar N}$}}\psfrag{S}{{\small $S$}}\psfrag{B}{{\small ${\bar S}$}}
\psfig{file=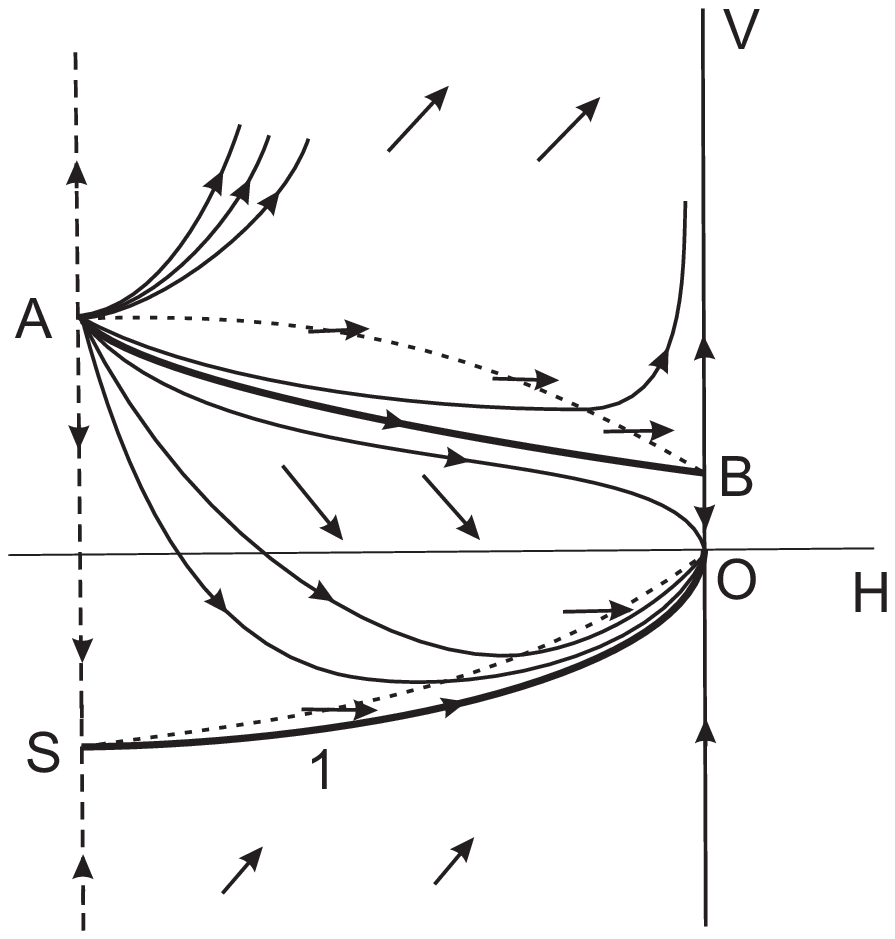,height=80mm,width=60mm} \centerline{\footnotesize
Figure \,8. The behavior of $w=w(h)$.}
\end{center}
\end{figure}

We proceed as in case 2). The graph of $w(h)$ coincides with the separatrix connecting the saddle
$S=(-\frac14, -\frac17)$  and the stable node $N$ at the origin. To verify this, we use Corollary 1 to obtain
$w=-\frac17+\frac{9}{28}s+\frac{99}{320}s^2+... $, hence $w(-\frac14)=-\frac17$, $w'(-\frac14)=\frac{9}{28}$ and
$w''(-\frac14)>0$. The horizontal isocline in the strip $-\frac14<h<0$ is composed of two parts of the left branch of the hyperbola
$\frac{35}{4}w^2-(2h+1)w+h=0$ going through the singular points $S$, $N$, the saddle $(0,\frac{4}{35})$ and the unstable node
$(-\frac14,\frac15)$. Its lower part $w=\psi(h)$, connecting $S$ and $N$ is increasing and $\psi'(-\frac14)=\frac37>w'(-\frac14)$.
Therefore $w(h)$ is increasing and convex in the whole interval and $w(h)<\psi(h)$, with the same proof as in Lemma 2. Then the
proof of 3) is the same as in Propositions 5 and 6.

\vspace{1ex}
\noindent
{\bf The proof of Theorem 3 (ii)} now follows from the picture shown in Figure 4. We take the case when $R_0(h)$ has 3 zeros and 3 extrema.
The other cases (less zeros or extrema, or double ones) would give the same or smaller result about the possible number of intersections.

\vspace{1ex}
\noindent
{\bf The proof of Theorem 3 (i)}. Since $R_\pm(h)=R_0(h) \pm (\beta_0+\beta_1 h)$, then $R_\pm'(h)=R_0'(h) \pm \beta_1$
 can have at most $3+4$ or $4+3$  zeros in $\Sigma$.  The proof is the same as above. Then $R_\pm(h)$  can have at most $4+5$ or $5+4$
total number of zeros.

The pictures drawn in Figures 4, 5 and 6 show respectively distributions $4+4$, $5+3$ and $5+4$. There, the cases when $R_0(h)$ has 3 zeros, 3 extrema and 3 inflection points are taken as giving the maximal cyclicity.
$\Box$


\section  {The small-amplitude limit cycles}

Here we use  Corollary 1 to estimate how many limit cycles system (\ref{E2}) can have in a vicinity of the centers $S$ and $S^*$.


\vspace{2ex}
\noindent
{\bf Proof of Theorem 4.}
With $h=s-\frac14$, we easily obtain that
$$M_2(h)=\bar{M}_2(s)=(\bar\alpha_0+\alpha_1s)I_0+(\bar\beta_0+\beta_1s)I_1+(\bar\gamma_0+\gamma_1s)I_2,\;\;
{\textstyle s\in (0,\frac14)}$$
where $\bar\alpha_0=\alpha_0-\frac14\alpha_1$ etc. By using the expansions from Corollary 1, one obtains by direct calculations
$$\bar{M}_2(s)=c_1\sum_{k=1}^\infty m_ks^k$$
where
$$\begin{array}{l}
m_1=\bar\alpha_0+\bar\beta_0+\bar\gamma_0,\\[2mm]
m_2=\frac38\bar\alpha_0+\alpha_1+\beta_1-\frac18\bar\gamma_0+\gamma_1,\\[2mm]
m_3=\frac{35}{64}\bar\alpha_0+\frac38\alpha_1-\frac{5}{64}\bar\gamma_0-\frac18\gamma_1,\\[2mm]
m_4=\frac{1155}{1024}\bar\alpha_0+\frac{35}{64}\alpha_1-\frac{105}{1024}\bar\gamma_0-\frac{5}{64}\gamma_1,\\[2mm]
m_5=\frac{45045}{16384}\bar\alpha_0+\frac{1155}{1024}\alpha_1-\frac{3003}{16384}\bar\gamma_0-\frac{105}{1024}\gamma_1.\\[2mm]
%
%
\end{array}$$
Assuming that $\gamma_1\neq 0$ which is the general situation (elsewhere one can produce less limit cycles), we first solve the equations
$m_3=m_4=m_5=0$ to obtain
$$\textstyle  \bar\alpha_0=-\frac{272}{539}\gamma_1,\;\;\alpha_1=\frac{1193}{539}\gamma_1,\;\;\bar\gamma_0=\frac{2960}{539}\gamma_1.$$
Then we solve equations $m_1=m_2=0$ to obtain
$${\textstyle \bar\beta_0=-\frac{384}{77}\gamma_1,\;\; \beta_1=-\frac{180}{77}\gamma_1,\;\;}
 \bar{M}_2(s)=c_1\gamma_1\sum_{k=6}^\infty d_ks^k\equiv  c_1\gamma_1R_6(s).$$
It is not difficult to see that $d_k>0$ if $k\geq 6$. Indeed, using the notations $a_k$ and $b_k$ from Corollary 1 and $m_k=\gamma_1d_k$
from above, we obtain for $k\geq 2$
$$\begin{array}{l}
c_1\gamma_1d_{k+1}=\bar\alpha_0a_{k+1}+\alpha_1a_k+\bar\gamma_0c_{k+1}+\gamma_1c_k\\[2mm]
=\frac{c_1\gamma_1}{539}\left[-272\frac{(4k-1)!!}{4^{k}k!(k+1)!}+1193\frac{(4k-5)!!}{4^{k-1}(k-1)!k!}-2960\frac{(4k-3)!!}{4^{k}k!(k+1)!}-539\frac{(4k-7)!!}{4^{k-1}(k-1)!k!}\right]\\[2mm]
= c_1\gamma_1\frac{15}{77}\frac{(4k-7)!!}{4^{k-2}k!(k+1)!}(k-2)(k-3)(k-4), \quad  k\geq 2.
\end{array}$$

In order to produce a perturbation with 5 small-amplitude limit cycles near $C$, we can keep $\gamma_1$ unchanged and move slightly the other
parameters of the present perturbation. This can be done by solving as before the above system where $m_k=0$, $k\leq 5$ are replaced
with $\hat{m}_k=\delta_k\gamma_1$, assuming that $\delta_k$ satisfy $\delta_1<0$, $\delta_k\delta_{k+1}<0$,  $|\delta_k|<\!\!<|\delta_{k+1}|$,
$|\delta_5|<\!\!<1$. Denote by $\hat{R}_6(s)>0$ the corresponding sum obtained with coefficients $\hat{d}_k$, $k\geq 6$ which slightly will differ from $d_k$.
Then one obtains a perturbation with
\begin{equation}\label{E22}
\bar{M}_2(s)=c_1\gamma_1[\delta_1s+\delta_2s^2+\delta_3s^3+\delta_4s^4+\delta_5s^5+\hat{R}_6(s)]
\end{equation}
 having 5 different small positive zeros.

The same perturbation will produce no limit cycles (neither small nor whatever) around the other center $C^*$ because by changing the signs of
 $\bar\beta_0$ and $\beta_1$, we will find that the corresponding constants to the left become $m_1^*=\frac{768}{77}\gamma_1$ and
 $m_2^*=\frac{360}{77}\gamma_1$ and therefore
$$\textstyle M^*_2(h)=\bar{M}^*_2(s)=c_1\gamma_1\left[\frac{768}{77}s+\frac{360}{77}s^2+R_6(s)\right]$$
will be positive and quite far from zero.

 Thus, the distributions $(5,0)$ and $(0,5)$ are only possible for  this case.  Below we are going to consider the cases when $\bar\beta_0$
or $\beta_1$ or both are zero.

Take first  $\bar\beta_0=\beta_1=0$. Solving the first three equations above with $m_1=m_2=m_3=0$ we obtain
$$\textstyle \bar\alpha_0=\frac{5}{28}\gamma_1,\; \;\alpha_1=-\frac{5}{56}\gamma_1, \;\;\bar\gamma_0=-\frac{5}{28}\gamma_1,\;\;
m_4=\frac{95}{1024}\gamma_1,\;\;\bar{M}_2(s)=c_1\gamma_1\frac{95}{1024}s^4[1+O(s)].$$
Then just as above we keep $\gamma_1$ the same and change slightly other parameters to obtain a system with
$\bar{M}_2(s)=c_1\gamma_1[\delta_1s+\delta_2s^2+\delta_3s^3+\hat{d}_4s^4(1+O(s))]$
having 3 small-amplitude limit cycles around $S$.  Since $M_2(h)=M_2^*(h)$, there are 3 limit cycles near $C^*$ too.

Let $\bar\beta_0=0$ now. Solving the first 4 equations with $m_1=m_2=m_3=m_4=0$ one obtains
$${\textstyle \bar\alpha_0=-\frac{16}{49}\gamma_1, \;\alpha_1=\frac{43}{49}\gamma_1, \;\bar\gamma_0=\frac{16}{49}\gamma_1,\;
 \beta_1=-\frac{12}{7}\gamma_1,\;}  \bar{M}_2(s)=c_1\gamma_1\sum_{k=5}^\infty d_ks^k=c_1\gamma_1 R_5(s).$$
Similarly to previous cases, one can calculate $d_k$  which  are all negative for $k\geq 5$ since
$$\textstyle d_{k+1}=-\frac37\frac{(4k-7)!!}{4^{k-2}k!(k+1)!!}k(k-2)(k-3),\quad k=2,3,...., d_5=-\frac{9}{128}.$$
As above, one can produce a secondary small perturbation with 4 small-amplitude limit cycles near $C$ and  $\bar{M}_2(s)$
taking the form
$$\bar{M}_2(s)=c_1\gamma_1[\delta_1s+\delta_2s^2+\delta_3s^3+\delta_4s^4+\hat{R}_5(s)]$$
where $\delta_k$ are chosen as in (\ref{E22}).
About $C^*$, by rotating the centers, one obtains a perturbation with $\bar{M}^*_2(s)=c_1\gamma_1\left[\frac{24}{7}s^2+R_5(s)\right].$
Clearly, the secondary small perturbation then leads to
$$\textstyle \bar{M}^*_2(s)=c_1\gamma_1[\delta_1s+(\frac{24}{7}+\delta_2)s^2+\delta_3s^3+\delta_4s^4+\hat{R}_5(s)], \;\;s\in(0,\frac14) $$
and $\bar{M}^*_2(s)$ will have a small positive zero close  to $-\frac{7}{24}\delta_1$. All the remaining terms are much smaller than the first two
in order to have another small zeros. Indeed, $R_5(s)$ (and therefore  $R_5(s)$) satisfies for $s\leq \frac18$ the inequality $|R_5(s)|\leq 2|d_5|s^5$.
This is because $|d_{k+1}|\leq 4|d_k|$.

And finally, if $\bar\beta_0\neq 0$, only one of the nests can have small-amplitude limit cycles even if $\beta_1=0$.  $\Box$


\vspace{2ex}
\section{Appendix}

Here we collect some identities we used in our proofs above. All they can be verified by direct calculations.
Let us recall that $H=\frac12y^2-\frac12x^2+\frac14 x^4$.

\vspace{1ex}
(I-1)  $ y^3dx=(\frac37x^2y+\frac{12}{7}yH)dx+\frac17dxy^3-\frac37xydH$


\vspace{1ex}
(II-1) $y^2dx=d(2xH+\frac13x^3-\frac{1}{10}x^5)-2xdH$

\vspace{1ex}
(II-2) $xy^2dx=d(x^2H+\frac14x^4-\frac{1}{12}x^6)-x^2dH$

\vspace{1ex}
(II-3) $x^3ydx=xydx-\frac13dy^3+ydH$

\vspace{1ex}
(II-4) $x^4ydx=(\frac87x^2y+\frac47yH)dx-\frac27dxy^3+\frac67xydH$

\vspace{1ex}
(II-5) $xy^3dx=\frac38(1+4H)xydx+\frac18d(x^2-1)y^3+\frac38(1-x^2)ydH$

\vspace{1ex}
(II-6) $x^2y^3dx=(\frac{4}{21}yH+\frac{8}{21}x^2y+\frac{4}{3}x^2yH)dx+d(\frac19x^3-
\frac{2}{21}x)y^3+(\frac27xy-\frac13x^3y)dH$

  \vspace{1ex}
(III-1) $x^2y^2dx=d(\frac23x^3H+\frac15x^5-\frac{1}{14}x^7)-\frac23x^3dH$

\vspace{1ex}
(III-2) $x^3y^2dx=d(\frac12x^4H+\frac16x^6-\frac{1}{16}x^8)-\frac12x^4dH$

\vspace{1ex}
(III-3)  $y^4dx=d(4xH^2+\frac43x^3H-\frac25x^5H+\frac15x^5-\frac17x^7+\frac{1}{36}x^9)$

\vspace{1ex}$\hspace{12mm}-(8xH+\frac43x^3-\frac25x^5)dH$

\vspace{1ex}
(III-4) $xy^4dx=d(2x^2H^2+x^4H-\frac13x^6H+\frac16x^6-\frac18x^8+\frac{1}{40}x^{10})$

\vspace{1ex}$\hspace{12mm}-(4x^2H+x^4-\frac13x^6)dH$

\vspace{1ex}
(III-5) $x^5ydx=(\frac54+H)xydx-d(\frac{5}{12}+\frac14x^2)y^3+(\frac54+\frac34x^2)ydH$

\vspace{1ex}
(IV-1) $y^4dq_1 \sim -[8\int q_1(x^3-x)dx+4y^2q_1]dH$

\vspace{1ex}
(IV-2) $q_1y^4dq_1 \sim -[4\int q_1^2(x^3-x)dx+2y^2q_1^2]dH$

\vspace{1ex}
(IV-3) $y^6dq_1 \sim [48\int q_1(x^3-x)(\frac14x^4-\frac12x^2)dx-48H\int q_1(x^3-x)dx-6y^4 q_1]dH$

\vspace{3ex}
\noindent
{\bf Acknowledgments.} Part of this work has been done while the first author visited Peking University during
August of 2018. He is very much grateful for excellent hospitality. The authors thank Lubomir Gavrilov  for his
useful comments.

\end{document}